\documentclass[11pt]{article}
\usepackage[T1]{fontenc}
\usepackage{amsthm,amsmath,amssymb,mathrsfs,extarrows,mathtools}
\usepackage{graphicx}
\usepackage{subfig}
\usepackage{float}
\usepackage{epstopdf}
\usepackage{graphics} 
\usepackage{epsfig} 
\usepackage{tikz}
\usepackage{indentfirst}
\usepackage{paralist}
\usepackage{booktabs}
\usepackage{multirow} 
\usepackage{multicol} 
\usepackage{enumerate}
\usepackage{enumitem}
\usepackage{mdwlist}
\usepackage{algorithm,algorithmicx}
\usepackage{hyphenat}          
\usepackage{microtype} 
\usepackage[colorlinks,linkcolor=blue,anchorcolor=blue,citecolor=red]{hyperref}
\makeatother
\newtheorem{theorem}{Theorem}[section]

\newtheorem{lemma}[theorem]{Lemma}

\newtheorem{remark}[theorem]{Remark}

\numberwithin{equation}{section}

\makeatletter

\newcommand{\Rmnum}[1]{\expandafter\@slowromancap\romannumeral #1@}
\makeatother

\newcommand{\me}{\mathrm{e}}
\newcommand{\mi}{\mathrm{i}}

\newcommand{\mr}{\mathbb{R}^2}

\newcommand{\fp}{\mathfrak{p}}
\newcommand{\sg}{\boldsymbol{\sigma}}
\newcommand{\fs}{\mathfrak{s}}

\newcommand{\od}{\overline{D}}
\newcommand{\oo}{\mathcal{O}}

\newcommand{\kp}{\kappa_\mathfrak{p}}
\newcommand{\ks}{\kappa_\mathfrak{s}}

\newcommand{\bx}{\mathbf{x}}
\newcommand{\bxs}{\mathbf{x_0}}

\newcommand{\bal}{\boldsymbol{\alpha}}

\newcommand{\grad}{\textrm{grad}}
\renewcommand{\div}{\textrm{div}}
\newcommand{\bu}{\boldsymbol{u}}

\newcommand{\bv}{\boldsymbol{v}}

\newcommand{\ls}{\xi_\fs}
\newcommand{\lp}{\xi_\fp}
\newcommand{\kpz}{\kappa_{\mathfrak{p},z}}
\newcommand{\ksz}{\kappa_{\mathfrak{s},z}}

\newcommand{\face}{\Gamma_0}
\newcommand{\lowsp}{\mathbb{R}^2_-}
\newcommand{\upsp}{\mathbb{R}^2_+}
\newcommand{\by}{\boldsymbol{y}}

\newcommand{\hx}{\hat{x}}

\newcommand{\uinc}{\boldsymbol{u^i} }
\newcommand{\ubd}{\boldsymbol{u^{bk}} }

\newcommand{\usc}{\boldsymbol {u_1^{\rm sc}}}
\newcommand{\uscp}{\boldsymbol {u_2^{\rm sc}}}
\newcommand{\uscd}{\boldsymbol {u^{\rm sc}_D}}
\newcommand{\uscn}{\boldsymbol {u^{\rm sc}_N}}
\newcommand{\es}{e_\fs}
\newcommand{\ep}{e_\fp}

\begin{document}
	\begin{titlepage}
		\title{A fast solver for many-particle elastic scattering in layered media}
		\author{Jinrui Zhang\thanks{Department of Mathematics, The Hong Kong University of Science and Technology,	Clear Water Bay, Kowloon, Hong Kong Special Administrative Region of China. Email: {\tt jinruizhang@ust.hk}.}\;,
        Yixiao He\thanks{School of Mathematical Sciences, Zhejiang University, Hangzhou, Zhejiang 310027, China. Email: {\tt \text{hyx\_math@zju.edu.cn}}.}\;,
        Jun Lai\thanks{School of Mathematical Sciences, Zhejiang University, Hangzhou, Zhejiang 310027, China, and Center for
Interdisciplinary Applied Mathematics, Zhejiang University, Hangzhou, Zhejiang, 310027, China. Email: {\tt laijun6@zju.edu.cn}.}
			}
		\date{}
	\end{titlepage}

		
	\maketitle
	\begin{abstract}
		This paper proposes a fast solver for time-harmonic elastic scattering by multiple particles embedded in layered media, with either Dirichlet or Neumann boundary conditions imposed on the particle surfaces.  Such problems arise in many important applications, including composite material optimization, nondestructive testing, and subsurface imaging. They are computationally challenging because of strong multiple scattering interactions among the particles and the layered interface. The proposed method represents the layered medium contribution through Sommerfeld integrals and couples this representation with a well-posed boundary integral formulation for the particle scattering problem. High-order integral equation discretization and scattering matrix are used to handle particles of general shape, while multiple scattering theory provides an efficient description of particle interactions. To reduce the cost for large particle systems, the resulting multiple scattering computation is further accelerated by the fast multipole method. The main formulation is developed in both two and three dimensions. Numerical experiments for rigid and traction-free particles validate the accuracy of the formulation, and demonstrate its flexibility in both direct scattering simulations and inverse scattering applications.
		\vspace{1em}
		
	\end{abstract}
	
	
	\section{Introduction}
	
	Wave scattering by many particles embedded in a layered medium arises in a broad range of scientific and engineering applications, including geophysical exploration, earthquake modeling, nondestructive testing,  and the design and optimization of composite materials \cite{achenbach2012wave,shearer2019introduction}. In these settings, the particles may model inclusions, cavities, defects, fibers, or small-scale heterogeneities, while the planar interface represents a sharp change in material properties. Accurate simulation of such configurations is important not only for forward prediction, but also for inverse problems in which the location, shape, or material character of buried objects must be inferred from measured wave fields.
	
	The computation becomes particularly challenging when the particle size is comparable to the incident wavelength and the particle distribution is dense. In this regime, multiple scattering among particles is strong, and effective-medium or homogenization approximations may no longer capture the relevant wave interactions \cite{parnell2010effective}. Instead, one must solve the full elastic scattering problem repeatedly, for instance during optimization, uncertainty quantification, or iterative inversion. Direct discretization of the full computational domain is prohibitively expensive for large particle systems because it must resolve both the local particle geometry and the long-range wave interactions mediated by the layered background. Thus, fast solvers that retain the accuracy of boundary integral formulations while scaling efficiently with the number of particles are essential.
	
	Substantial progress has been made on fast algorithms for acoustic, electromagnetic, and elastic multiple scattering in homogeneous media. Related fast multiple particle solvers have been developed for acoustic scattering \cite{laiFastSolverMultiparticle2014}, elastic scattering in two- and three-dimensional free space \cite{laiFrameworkSimulationMultiple2019,laiFastInverseElastic2022}, and electromagnetic scattering \cite{gimbutasFastMultiparticleScattering2013}.  For elastic waves in half-space or rough-surface configurations, several analytical and numerical studies have also considered Dirichlet or Neumann boundary conditions imposed on the interface \cite{huDirectInverseElastic2018,huElasticScatteringRough2020}. Existing numerical approaches mainly rely on explicit layered Green's tensors, whose repeated evaluation is often expensive \cite{touheiAnalysisScatteringWaves2003,liuThreedimensionalIndirectBoundary2019}, or, more recently, on windowed Green's function \cite{brunoWindowedGreenFunction2021,yin2025efficient}. However, efficient solvers for large-scale many-particle elastic scattering in layered media remain largely unexplored.  For acoustic and electromagnetic half-space problems with homogeneous Dirichlet or Neumann boundary conditions on a flat interface, reflection principles can often be employed to reduce the problem to an equivalent free-space formulation. In contrast, the Navier equation subject to homogeneous boundary conditions does not admit such a straightforward reflection construction. Furthermore, for penetrable interface problems, the compressional and shear wave components are coupled through the continuity conditions for displacement and traction, making the analysis substantially more involved than in the acoustic or electromagnetic setting \cite{laiFastSolverMultiparticle2014}. 
    	
	In this paper, we develop a fast solver for time-harmonic elastic scattering by a large number of rigid or traction-free particles embedded in the lower half-space of a two-layered medium, with penetrable transmission conditions imposed on the flat interface. The response of the layered background is represented by Sommerfeld integrals, which account for the reflection and transmission of both compressional and shear waves at the interface. Scattering by individual particles, including particles of non-circular shape, is resolved through a well-posed boundary integral formulation and high-order discretization. Scattering matrix theory is then used to express the incoming and outgoing multipole coefficients of each particle, thereby separating local particle scattering from global particle-to-particle interactions. The free-space part of the many-body interaction is accelerated by the fast multipole method, and analytical conversion formulas connect the multipole representation with the Sommerfeld integral representation required by the layered medium. This framework yields a solver that is accurate for complex particle geometries and efficient for large-scale configurations.
	
	The main contribution of this work is a unified multiple-scattering framework that incorporates the elastic layered Green's function, scattering matrices for general particles, and fast multipole acceleration. By explicitly accounting for field and traction continuity across the infinite interface, the proposed approach extends fast multiple particle scattering techniques from homogeneous elastic media to layered elastic media. Although we focus on a two-layer geometry, the construction is compatible with extensions to more general layered backgrounds by the response matrix for stratified materials \cite{koScatteringMatrixMethod1988}. Numerical examples demonstrate the accuracy and efficiency of the solver for direct scattering and illustrate its potential in inverse scattering applications, including selective focusing of buried elastic particles.
	
	An outline of the paper follows. In Section \ref{formulation}, we formulate the layered elastic scattering problem. The Sommerfeld integral in a layered medium is introduced in Section \ref{sommerfeld}. In Section \ref{scadisk}, we review the classical multiple scattering theory for circular particles and extend the scattering formalism to non-circular elastic particles. Section \ref{multi_part_layer} extends the fast multiple scattering method to layered media by developing analytical tools needed to connect the Sommerfeld integral formalism with multiple scattering theory. An  extension to three-dimensional setting is given in Section \ref{extension3D}. Numerical examples, including direct and inverse scattering tests, are provided in Section \ref{numer_exp} to illustrate the efficiency of the method, followed by concluding remarks in Section \ref{conc}.

	\section{Formulation of layered elastic scattering}\label{formulation}
\begin{figure}
		\centering
		\includegraphics[width =10cm]{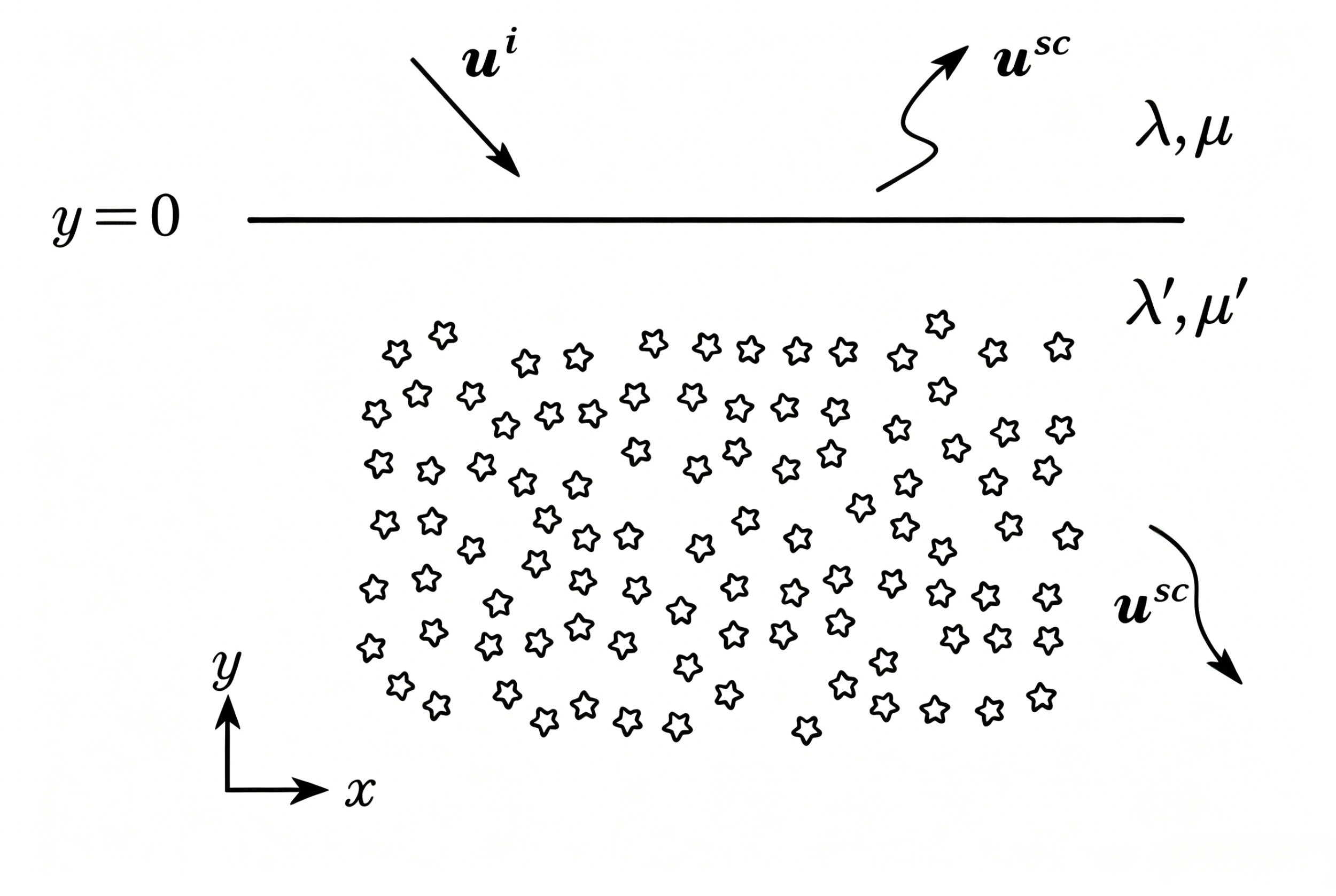}
		\caption{Geometry of the layered medium with many elastic particles buried in $\lowsp$.}\label{f1}
	\end{figure}

	For clarity, we begin with the two-dimensional formulation and defer the extension to three dimensions to Section \ref{extension3D}. For $\bx=(x_1,x_2)\in\mr$, let $\lowsp=\{(x_1,x_2)\in\mr\mid x_2<0\}$ and $\upsp=\{(x_1,x_2)\in\mr\mid x_2>0\}$ denote the lower and upper half-spaces, respectively. The interface separating the two layers is denoted by $\face=\{(x_1,x_2)\in\mr\mid x_2=0\}$. 
	Consider $M$ well-separated impenetrable particles fully embedded in the lower half-space $\lowsp$, denoted by $D=D_{1}\cup D_{2}\cup\cdots\cup D_{M}$, with Lipschitz boundary $\partial D$. See Figure \ref{f1} for an illustration.
	Let $\nu$ denote the unit exterior normal vector on $\Gamma$. When no ambiguity arises, the same symbol is also used for the upward unit normal vector on $\face$. 
	Let the particles be illuminated by a time-harmonic incident wave $\uinc$. The total displacement field $\boldsymbol u$ is the sum of the background field $\ubd$(including the reflected field and transmitted field in the absence of particles \cite{laizhang2026}), and the scattered field $\boldsymbol v$, namely $\boldsymbol u=\ubd+\boldsymbol v$, and satisfies the Navier equation
	\begin{equation*}
		\mu\Delta\boldsymbol{u}+(\lambda+\mu)\nabla\nabla\cdot\boldsymbol{u}
		+\rho\omega^2\boldsymbol{u}=0,\mbox{ in } \mathbb{R}^2\setminus \overline{D}, 	
	\end{equation*}
	where $\omega>0$ denotes the angular frequency, $\rho$ denotes the mass density, and $\lambda$ and $\mu$ are the piecewise constant Lam\'{e} parameters, which differ between the two layers and satisfy $\mu>0$ and $\lambda+\mu>0$. Without loss of generality, we assume that the elastic mass density $\rho$ is normalized to one throughout the medium.
	For a nonzero vector $x\in\mr$, we introduce the unit vector $\hx:=\bx/|\bx|$ and denote by $\hx^\bot$ the vector obtained by rotating $\hx$ counterclockwise by $\pi/2$.
	
	In addition to the usual differential operators $\grad$ and $\div$, we use the rotated gradient and divergence operators
	\[\grad^\bot \bu=\nabla^\bot \bu:=\left[-\frac{\partial u}{\partial x_2},\frac{\partial u}{\partial x_1}\right]^\top,\qquad \div^\bot\bu=\nabla^\bot\cdot\bu:=\frac{\partial u_2}{\partial x_1}-\frac{\partial u_1}{\partial x_2}.\]
	The linearized strain tensor is defined by 
	\[\varepsilon(\bu):=\frac{1}{2}(\nabla \bu+\nabla \bu^\top)\in\mathbb{C}^{2\times2},\]
	where $\nabla \bu$ and $\nabla \bu^\top$ denote the Jacobian matrix of $\bu(x)\in \mathbb{C}^2$ and its transpose, respectively. By Hooke's law, the strain tensor is related to the stress tensor through
	\begin{align*}
		\sigma(u)=\lambda(\nabla\cdot \bu)\mathbb{I}+2\mu\varepsilon(\bu)\in\mathbb{C}^{2\times2},
	\end{align*}
	where $\mathbb{I}$ is the $2\times 2$ identity matrix. The surface traction on $\Gamma$ is defined by
	\begin{eqnarray}\label{tracdef}
		T_\nu \bu: = \nu\cdot\sigma(\bu)=2\mu \nu\cdot \nabla\bu + \lambda \nu\nabla\cdot\bu - \mu\nu^\bot\nabla^\bot\cdot\bu.
	\end{eqnarray}
	
	
	It follows that the scattered field $\boldsymbol v$ satisfies the boundary value problem
	\begin{equation}\label{scatteredfield}
		\begin{cases}
			\mu\Delta\boldsymbol{v}+(\lambda+\mu)\nabla\nabla\cdot\boldsymbol{v}
			+\omega^2\boldsymbol{v}=0\quad &{\rm in}~
			\mr\setminus\od,\\
			\mathcal{B}(\boldsymbol{v})=-\mathcal{B}(\ubd)\quad &{\rm on}~\partial D,\\
			\left[\bv\right]=0,\left[T_\nu\bv\right]=0 &{\rm on}~\face,
		\end{cases}
	\end{equation}
	where $\left[\cdot\right]$ denotes the jump across the interface $\face$, and $\mathcal{B}$ represents one of the following particle boundary conditions:
	\begin{align*}
		\mbox{ rigid: }\mathcal{B}(\bv)&=\bv\quad{\rm on}~\Gamma,\\
		\mbox{ traction-free: }\mathcal{B}(\bv)&=T_\nu\bv\quad{\rm on}~\Gamma .
	\end{align*}
	
	By the Helmholtz decomposition, any solution $\bv$ of \eqref{scatteredfield} can be decomposed as
	\begin{align}\label{helmdec}
		\bv = \bv_\fp+\bv_\fs,\quad\boldsymbol v_{\mathfrak p}=-\frac{1}{\kappa_{\mathfrak p}^2}\nabla\nabla\cdot\boldsymbol v,\quad
		\boldsymbol v_{\mathfrak s}=-\frac{1}{\kappa_{\mathfrak s}^2}\nabla^\bot\nabla^\bot\cdot\boldsymbol v, 
	\end{align}
	where the compressional wavenumber $\kp$ and shear wavenumber $\ks$ are given by
	$
	\kappa_{\mathfrak p}={\omega}/{\sqrt{\lambda+2\mu}},
	\kappa_{\mathfrak s}={\omega}/{\sqrt{\mu}} 
	$.
	Here $\bv_\fp$ and $\bv_\fs$ are the compressional and shear components of $\boldsymbol v$, respectively. In addition, the scattered field is required to satisfy the Kupradze--Sommerfeld radiation condition
	\[
	\lim_{r\to\infty}\sqrt{r}(\partial_r\boldsymbol v_{\mathfrak p}-\mathrm{i}\kappa_{\mathfrak p}\boldsymbol v_{\mathfrak p})=0,\quad
	\lim_{r\to\infty}\sqrt{r}(\partial_r\boldsymbol v_{\mathfrak s}-\mathrm{i}\kappa_{\mathfrak s}\boldsymbol v_{\mathfrak s})=0,\quad r=|\bx|.
	\]
	
	To develop an efficient fast solver, we impose two additional assumptions. First, the particles $\{D_j\}$ are assumed to be well separated. More precisely, the distance between any two particles is taken to be at least 10\% of the particle size. This condition ensures that the multiple scattering representation remains accurate and stable. Second, we assume that only finitely many distinct particle shapes appear in the simulation. Under this assumption, the single-particle scattering matrices can be precomputed and reused, leading to a substantial reduction in the number of degrees of freedom required by the solver. The particles need not be symmetric and may be placed with arbitrary orientations. Both assumptions are natural in many applications involving repeated inclusions or defects in a layered background.
	\section{The Sommerfeld integral for elastic layered media}\label{sommerfeld}
	Wave propagation in layered media is a classical topic in acoustic, electromagnetic, and elastic scattering theory. A standard tool for such problems is the Sommerfeld representation, which expresses the field as a spectral integral in the transverse variable \cite{chew1999waves}. This representation is especially useful for planar interfaces because the transmission conditions can be enforced mode by mode in the Fourier domain. Throughout the following sections, the Lam\'e constants, compressional wavenumber, and shear wavenumber in the upper half-space $\upsp$ are denoted by $\lambda,\mu,\kp$, and $\ks$, respectively, while the corresponding quantities in the lower half-space $\lowsp$ are denoted by $\lambda',\mu',\kp'$, and $\ks'$, as shown in Figure \ref{f1}.
	
	\subsection{The free-space Green's function in Sommerfeld form}
	Let $\mathbf{x_0}=(x_0,y_0)$, $\mathbf x=(x,y)$ be the source and target points, respectively. The two-dimensional free-space elastic Green's function is given by \cite{hsiaoBoundaryIntegralEquations2008}:
	\begin{eqnarray}\label{greenfun}
		\Phi(\mathbf x,\mathbf{x_0}) = \frac{\mi\ks^2}{4\omega^2}H_0^{(1)}(\ks|\mathbf x-\mathbf{x_0}|)\mathbb{I}+\frac{\mi}{4 \omega^2}\nabla\nabla^\top \left[H_0^{(1)}(\ks|\mathbf x-\mathbf{x_0}|)-H_0^{(1)}(\kp|\mathbf x-\mathbf{x_0}|)\right],
	\end{eqnarray}
	where $\mathbb{I}$ is the $2\times 2$ identity matrix, and $H_0^{(1)}$ is the Hankel function of the first kind of order zero.
	Using the Fourier transform together with contour integration \cite{oneilEfficientRepresentationHalfspace2014}, one obtains
	\begin{align}\label{som1}
		\frac{\mi}{4}H_0^{(1)}(k|\mathbf{x}-\mathbf{x_0}|)=\frac1{4\pi}\int_{-\infty}^\infty\frac{\me^{-\sqrt{\xi^2-{k}^2}|y-y_0|}}{\sqrt{\xi^2-{k}^2}}\me^{\mi\xi(x-x_0)}d\xi.
	\end{align}
	The Sommerfeld integral in \eqref{som1} is conditionally convergent and is valid for $y\neq y_0$.
	Denote the spectral density associated with $\Phi(\mathbf x,\mathbf{x_0})$ by
	\[\widetilde{\Phi}(\xi,\mathbf x,\mathbf{x_0}) =\frac{1}{\omega^2} \left({\ks^2}\widetilde{g}_{\ks}\mathbb{I}+\nabla\nabla^\top(\widetilde{g}_{\ks}-\widetilde{g}_{\kp}) \right),\]
	where \[\widetilde{g}_k(\xi,\mathbf x,\mathbf{x_0}) =\frac{1}{4\pi} \frac{\me^{-\sqrt{\xi^2-{k}^2}|y-y_0|}}{\sqrt{\xi^2-{k}^2}}\me^{\mi\xi(x-x_0)}.\]
	Then, the free-space Green's tensor \eqref{greenfun} can be expressed as the Sommerfeld integral
	\[\Phi(\mathbf x,\mathbf{x_0}) =\int_{-\infty}^\infty\widetilde{\Phi}(\xi,\mathbf x,\mathbf{x_0})d\xi.\]
	\subsection{Elastic point source in a two-layer structure}
	Consider an elastic point source located at $\mathbf{x_0}=(x_0,y_0)$ in $\upsp$ with polarization direction $\bal = (\alpha_1,\alpha_2)^\top$. The incident field is
	\[\uinc = \Phi(\mathbf x,\mathbf{x_0})\bal =\int_{-\infty}^\infty\widetilde{\Phi}(\xi,\mathbf x,\mathbf{x_0})\bal d\xi.\]
	In the Sommerfeld formulation, the upward scattered field $\bu^s_1$ in the upper half-space $\upsp$ is represented as
	\begin{align}\label{us1}
		\bu^s_1 = \int_{-\infty}^\infty\widetilde{\Phi}^1(\xi,\mathbf x,\mathbf{x_0})\sg_1(\xi)d\xi,
	\end{align}
	where $\sg_1(\xi) = (\sigma_{11},\sigma_{12})^\top$ is an unknown upward spectral density associated with the interface $\face$, and
	\[\widetilde{\Phi}^1(\xi,\mathbf x,\mathbf{x_0}) =\frac{1}{\omega^2} \left({\ks^2}\widetilde{g}_{\ks}^1\mathbb{I}+\nabla\nabla^\top(\widetilde{g}_{\ks}^1-\widetilde{g}_{\kp}^1)\right),\quad  \widetilde{g}_k^1(\xi,\mathbf x,\mathbf{x_0}) =\frac{1}{4\pi} \frac{\me^{-\sqrt{\xi^2-{k}^2}y}}{\sqrt{\xi^2-{k}^2}}\me^{\mi\xi(x-x_0)} .\]
	It is straightforward to verify that $\bu^s_1$ satisfies the Navier equation in $\upsp$.
	Similarly, the downward scattered field $\bu_2^s$ in the lower half-space $\lowsp$ is written as
	\begin{align}\label{us2}
		\bu^s_2 = \int_{-\infty}^\infty\widetilde{\Phi}^2(\xi,\mathbf x,\mathbf{x_0})\sg_2(\xi)d\xi,
	\end{align}
	where $\sg_2(\xi) = (\sigma_{21},\sigma_{22})^\top$ is the corresponding unknown downward spectral density on $\face$, and
	\begin{align}\label{lkel}
		\widetilde{\Phi}^2(\xi,\mathbf x,\mathbf{x_0}) =\frac{1}{\omega^2} \left({{\ks'}^2}\widetilde{g}_{\ks'}^2\mathbb{I}+\nabla\nabla^\top(\widetilde{g}_{\ks'}^2-\widetilde{g}_{\kp'}^2) \right),\quad
		\widetilde{g}_k^2(\xi,\mathbf x,\mathbf{x_0}) =\frac{1}{4\pi} \frac{\me^{\sqrt{\xi^2-{k}^2}y}}{\sqrt{\xi^2-{k}^2}}\me^{\mi\xi(x-x_0)}.
	\end{align}
	\begin{remark}
		The signs of the factors $\me^{\pm\sqrt{\xi^2-{k}^2}y}$ in $\widetilde{g}_k^1$ and $\widetilde{g}_k^2$ ensure that evanescent modes, corresponding to $|\xi|>|k|$, decay away from the interface. 
	\end{remark}
	\begin{remark}
	The two unknown densities $\sg_1$ and $\sg_2$ admit two equivalent interpretations. They may be viewed as Fourier-domain spectral densities in a representation satisfying the Navier equation in each half-space. Equivalently, from the perspective of potential theory \cite{CotKress1983}, they can be interpreted as the Fourier transforms of source densities for two single-layer potentials supported on the interface $\face$.
		\end{remark}
	In the absence of elastic particles, the unknown functions $\sg_1$ and $\sg_2$ are determined by enforcing continuity of displacement and traction across the interface:
	\[\uinc+\bu_1^s|_{y=0}=\bu_2^s|_{y=0},\quad T_y(\uinc+\bu_1^s)|_{y=0}=T'_y(\bu_2^s)|_{y=0},\]
	where $T_y$, defined by \eqref{tracdef}, is the traction operator on $\face$ with $\nu=(0,1)$. In this case, it reduces to
	\begin{eqnarray*}
		T_y \bu=\begin{pmatrix}
			\mu\frac{ \partial u_1}{\partial y}+\mu\frac{ \partial u_2}{\partial x},
			(\lambda+2\mu)\frac{ \partial u_2}{\partial y}+\lambda\frac{ \partial u_1}{\partial x}
		\end{pmatrix}^\top, \quad \bu = (u_1,u_2)^\top.
	\end{eqnarray*}
	The operator $T'_y$ is obtained from $T_y$ by replacing $\lambda,\mu$ with $\lambda',\mu'$.
	
	Because the Sommerfeld representation diagonalizes the tangential variable, the interface conditions can be imposed independently for each Fourier mode $\xi$. A direct calculation gives the following linear system for $\sg=(\sigma_{11},\sigma_{12},\sigma_{21},\sigma_{22})^\top$:
	\begin{align}\label{ae}
		A_\xi\sg=\mathbf{b}_\xi,
	\end{align}
	$A_\xi=$
	{\small
		\begin{align}\label{aa}
			\begin{pmatrix}
				\frac{\ks^2}{\ls}-\frac{\xi^2}{\ls}+\frac{\xi^2}{\lp}&0&
				\frac{\xi^2}{\ls'}-\frac{\ks'^2}{\ls'}-\frac{\xi^2}{\lp'}&0\\ 0&\frac{\ks^2}{\ls}+\ls-\lp&\lp'-\frac{\ks'^2}{\ls'}-\ls'&0\\
				-\mu \ks^2&\mi\mu\xi\left(\frac{\ks^2}{\ls}+2\ls-2\lp\right)&-\mu' {\ks'}^2&\mi\mu'\xi\left(2\lp'-2\ls'-\frac{{\ks'}^2}{\ls'}\right)\\
				\mi\xi\left(2\mu (\ls-\lp)+\lambda\frac{\kp^2}{\lp}\right)&-(\lambda+2\mu){\kp}^2&\mi\xi\left(2\mu' (\lp'-\ls')-\lambda'\frac{{\kp'}^2}{\lp'}\right)&-(\lambda'+2\mu'){\kp'}^2
			\end{pmatrix},
		\end{align}
		\begin{align*}
			\mathbf{b}_\xi = -\begin{pmatrix}
				\left(-\ls\es+\frac{\xi^2\ep}{\lp}\right)\alpha_1+\mi\xi(\es-\ep)\alpha_2\vspace{1em}\\
				\mi\xi(\es-\ep)\alpha_1+\left( \frac{\xi^2\es}{\ls}-{\ep}{\lp}\right)\alpha_2\vspace{1em}\\
				\mu \left(\ks^2 \es-2{\xi^2\es}+2{\xi^2\ep}\right)
				\alpha_1
				+\mi\mu\xi\left(\ks^2 \frac{\es}{\ls}+2{\es}{\ls}-2{\ep}{\lp}\right)
				\alpha_2\vspace{1em}\\
				\mi\xi\left(2\mu (\ls \es-\lp \ep)+\lambda\kp^2\frac{\ep}{\lp}\right)\alpha_1
				+\left(2\mu\xi^2\left(\es-\ep\right)+(\lambda+2\mu)\kp^2\ep\right)\alpha_2
			\end{pmatrix},
	\end{align*}}
	where
	$\ls = \sqrt{\xi^2-{\ks}^2},\ls'= \sqrt{\xi^2-{\ks'}^2},
	\lp = \sqrt{\xi^2-{\kp}^2},\lp'= \sqrt{\xi^2-{\kp'}^2},\es = \me^{-\ls y_0},\ep = \me^{-\lp y_0}$.
	For the scattering problems considered in this paper, the Sommerfeld integrals must be coupled with a representation of the field generated by multiple particles buried in the lower half-space $\lowsp$. Before introducing the coupled layered-medium system, we first present in the next section the scattering formalism for a finite collection of elastic particles in a homogeneous background.
		\section{Wave scattering for multiple elastic particles}\label{scadisk}
	In this section, we introduce a fast numerical method for elastic scattering by multiple particles randomly distributed in a homogeneous and isotropic elastic background. The method is based on classical Mie theory and multipole translation formulas, and it provides the basic multiple-scattering framework that will later be coupled with the layered-medium Sommerfeld representation \cite{laiFrameworkSimulationMultiple2019}.
	\subsection{Scattering of a single disk}\label{singds}
	Consider an elastic disk embedded in a homogeneous medium with Lam\'e constants $\lambda$ and $\mu$ and angular frequency $\omega$. The corresponding compressional and shear wavenumbers are $\kp$ and $\ks$, respectively. Let the disk $S_0$ be centered at the origin with radius $R$. 
	For a given point $\bx=(x_1,x_2)$, denote its polar coordinates by $(r,\theta)$. Let $J_n$ and $H_n^{(1)}$ be the Bessel and Hankel functions of the first kind of order $n$, respectively.
	Define the scalar functions
	\begin{eqnarray*}
		u_{n}^\kappa (\bx) = J_n(\kappa r)\me^{\mi n\theta },v_{n}^\kappa (\bx) = H_n^{(1)}(\kappa r)\me^{\mi n\theta },
	\end{eqnarray*}
	which are referred to as \textit{cylindrical wave functions}. They satisfy the two-dimensional Helmholtz equation, with $v^{\kappa}_{n}(\bx)$ singular at the origin.
	
	According to the Helmholtz decomposition \eqref{helmdec}, the incoming field near the disk $S_0$ can be expanded as
	\begin{eqnarray*}
		\uinc(\bx) = \sum_{n=-\infty}^{\infty}a_{n}\nabla u_{n}^{\kp} (\bx)+b_{n}\nabla^\bot u_{n}^{\ks} (\bx),
	\end{eqnarray*}
	where $\{a_{n},b_{n}\}$ are the \textit{incoming expansion coefficients} of $\uinc$ with respect to $S_0$.
	
	The corresponding scattered field $\bv$ in the exterior of $S_0$ has the outgoing expansion
	\begin{eqnarray}\label{scatspan}
		\bv =\sum_{n=-\infty}^{\infty}\alpha_{n}\nabla v_{n}^{\kp} (\bx)+\beta_{n}\nabla^\bot v_{n}^{\ks} (\bx),
	\end{eqnarray} 
	where $\{\alpha_{n},\beta_{n}\}$ are the \textit{outgoing expansion coefficients}. The linear map from the incoming coefficients $\{a_{n},b_{n}\}$ to the outgoing coefficients $\{\alpha_{n},\beta_{n}\}$, $n\in\mathbb{Z}$, is called the \textit{scattering matrix} $\mathscr{S}$. For a disk, this matrix is obtained by enforcing the boundary condition mode by mode on the disk boundary. Detailed derivations for rigid and traction-free disks can be found in \cite{laiFrameworkSimulationMultiple2019} and \cite{zhangSelectiveFocusingElastic}, respectively.
	\begin{remark}
		In practice, the infinite expansions are truncated to modes indexed by $\mathbb{Z}_p = \{-p,-p+1,\ldots,0,1,\ldots,p\}$. The choice of the truncation parameter $p$ has been studied numerically in \cite{antoineNumericalApproximationHighfrequency2008,antoineWideFrequencyBand2012}. For any coefficient sequence $\gamma=a,b,\alpha$, or $\beta$, we write $\vec{\gamma}\equiv(\gamma_{-p},\gamma_{-p+1},\ldots,\,\gamma_0,\gamma_1,\ldots,\gamma_p)^\top$.
	\end{remark}
	\subsection{Scattering of multiple disks}
	Suppose that $M$ well-separated, identical elastic disks of radius $R$ are distributed in a homogeneous medium. For each individual particle, the single-disk analysis described in Subsection \ref{singds} applies.
	Let $\vec{a}^m, \vec{b}^m$ denote the incoming coefficients and $\vec{\alpha}^m, \vec{\beta}^m$ the outgoing coefficients for the $m$-th particle. Then
	\begin{equation}\label{mup1}\begin{bmatrix}\vec{\alpha}^m\\ \vec{\beta}^m\end{bmatrix}=\mathscr{S}_p\begin{bmatrix}\vec{a}^m\\ \vec{b}^m\end{bmatrix}, m = 1,\cdots, M,
	\end{equation}
	where $\mathscr{S}_p$ denotes the truncated $(4p+2)\times(4p+2)$ scattering matrix acting on the truncated expansion.
	
	The principal difference between single-particle and multiple-particle scattering is that the incoming field experienced by each particle consists of two contributions: the applied incident field $\uinc$ and the fields scattered from all other particles. In particular, let us denote by $\mathscr{T}^{ml}$ the translation operator that maps the outgoing coefficients $\{\vec{\alpha}^l,\vec{\beta}^l\}$ of particle $l$ to the local incoming expansion $\{\vec{a}^m,\vec{b}^m\}$ centered at particle $m$. Then the incoming coefficients for the $m$-th particle satisfy
		\begin{equation}\label{mulp2}
			\begin{bmatrix}\vec{a}^m\\ \vec{b}^m\end{bmatrix}=\begin{bmatrix}\vec{a_0}^m\\ \vec{b_0}^m\end{bmatrix}+\sum_{l=1\atop l\neq m}^M\mathscr{T}^{ml}\begin{bmatrix}\vec{\alpha}^l\\ \vec{\beta}^l\end{bmatrix},
		\end{equation}
		where $\left\{\vec{a_0}^m,\vec{b_0}^m\right\}$ is the truncated local expansion of the incident wave $\uinc$ about particle $m$. The operator $\mathscr{T}^{ml}$ is the multipole-to-local (M2L) translation operator \cite{rokhlinRapidSolutionIntegral1990}. By Graf's addition theorem, the translation matrix has the form
		\[\mathscr{T}^{ml}=\begin{bmatrix}\{a_{ij}^{ml}\}_{i,j\in\mathbb{Z}_p}&0\\0&\{b_{ij}^{ml}\}_{i,j\in\mathbb{Z}_p}\end{bmatrix},\]
		where \[a_{ij}^{ml}=H_{i-j}^{(1)}(\kp|x_l-x_m|)\me^{-\mathrm{i}(i-j)(\theta_{ml}-\pi)},\quad b_{ij}^{ml}=H_{i-j}^{(1)}(\ks|x_l-x_m|)\me^{-\mathrm{i}(i-j)(\theta_{ml}-\pi)}.\]
	Combining \eqref{mup1} and \eqref{mulp2}, the incoming coefficients $\vec{a}^m,\vec{b}^m$ can be eliminated, yielding a linear system involving only the outgoing coefficients:
	\begin{align}\label{mulfor}
		\begin{bmatrix}\mathscr{S}_p^{-1}&-\mathscr{T}^{12}&\ldots&-\mathscr{T}^{1M}\\-\mathscr{T}^{21}&\mathscr{S}_p^{-1}&\ldots&-\mathscr{T}^{2M}\\\vdots&\vdots&\ddots&\vdots\\-\mathscr{T}^{M1}&-\mathscr{T}^{M2}&\ldots&\mathscr{S}_p^{-1}\end{bmatrix}
		\begin{bmatrix}\begin{bmatrix}\vec{\alpha}^1\\\vec{\beta}^1\end{bmatrix}\vspace{0.6em}\\
			\begin{bmatrix}\vec{\alpha}^2\\\vec{\beta}^2\end{bmatrix}\\
			\vdots\\
			\begin{bmatrix}\vec{\alpha}^M\\\vec{\beta}^M\end{bmatrix}\end{bmatrix}
		= 
		\begin{bmatrix}\begin{bmatrix}\vec{a_0}^1\\\vec{b_0}^1\end{bmatrix}\vspace{0.6em}\\
			\begin{bmatrix}\vec{a_0}^2\\\vec{b_0}^2\end{bmatrix}\\
			\vdots\\
			\begin{bmatrix}\vec{a_0}^M\\\vec{b_0}^M\end{bmatrix}\end{bmatrix}.
	\end{align}
	The system \eqref{mulfor} can be solved iteratively using GMRES \cite{saadGMRESGeneralizedMinimal1986}. Its conditioning has been studied for both dense media \cite{thierrySpectralConditionNumber2013} and dilute media \cite{antoineSpectralConditionNumber2013}. In practice, it is advantageous to precondition \eqref{mulfor} using a block-diagonal preconditioner whose diagonal blocks are the single-particle scattering matrix $\mathscr{S}_p$.
	The resulting system is
	\begin{align}\label{mul2}
		\begin{bmatrix}I&-\mathscr{S}_p\mathscr{T}^{12}&\ldots&-\mathscr{S}_p\mathscr{T}^{1M}\\-\mathscr{S}_p\mathscr{T}^{21}&I&\ldots&-\mathscr{S}_p\mathscr{T}^{2M}\\\vdots&\vdots&\ddots&\vdots\\-\mathscr{S}_p\mathscr{T}^{M1}&-\mathscr{S}_p\mathscr{T}^{M2}&\ldots&I\end{bmatrix}
		\begin{bmatrix}\begin{bmatrix}\vec{\alpha}^1\\\vec{\beta}^1\end{bmatrix}\vspace{0.6em}\\
			\begin{bmatrix}\vec{\alpha}^2\\\vec{\beta}^2\end{bmatrix}\\
			\vdots\\
			\begin{bmatrix}\vec{\alpha}^M\\\vec{\beta}^M\end{bmatrix}\end{bmatrix}
		= 
		\begin{bmatrix}\mathscr{S}_p\begin{bmatrix}\vec{a_0}^1\\\vec{b_0}^1\end{bmatrix}\vspace{0.6em}\\
			\mathscr{S}_p\begin{bmatrix}\vec{a_0}^2\\\vec{b_0}^2\end{bmatrix}\\
			\vdots\\
			\mathscr{S}_p\begin{bmatrix}\vec{a_0}^M\\\vec{b_0}^M\end{bmatrix}\end{bmatrix}.
	\end{align}
	This preconditioned system is typically much better conditioned than the original system \eqref{mulfor} and therefore converges faster. 
	Furthermore, since each translation operator $\mathscr{T}^{ml}$ is dense, a direct matrix-vector product has computational complexity of order $\oo(M^2p^2)$. The fast multipole method (FMM) reduces this cost to $\oo(Mp^2)$ per iteration and substantially accelerates the computation \cite{carrierFastAdaptiveMultipole1988,laiimpleFastMultipole2025,gimbutasFastMultiparticleScattering2013}. The convergence of such multipole expansion methods has recently been analyzed in \cite{fitzpatrickConvergenceMultipoleExpansion2021}.
	

	\subsection{Scattering of arbitrarily shaped particles}\label{arbq}

\begin{figure}
		\centering
		\includegraphics[scale =0.3]{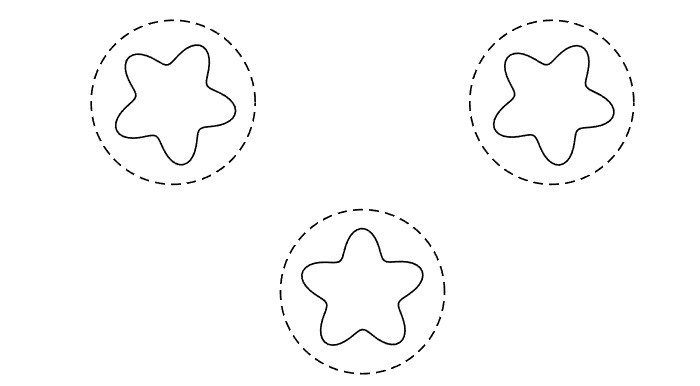}
		\caption{Elastic scattering by multiple obstacles with enclosing disks.}\label{f2}
	\end{figure}

	For arbitrarily shaped elastic particles, the disk-based multiple scattering theory described above cannot be applied directly. In this section, we extend the scattering matrix framework to noncircular particles. The construction follows our previous fast algorithm for many arbitrarily shaped elastic particles \cite{laiFrameworkSimulationMultiple2019}, with two modifications. First, we use the high-precision boundary integral solver in \cite{yaoRobustHighPrecision2024}, which requires a new formula that maps the boundary density to outgoing expansion coefficients. Second, we also include the traction-free case. We assume that the particles are well separated in the sense that each particle is contained in an enclosing disk and that these disks do not overlap. An illustration is shown in Figure \ref{f2}.
	For each particle, we construct a scattering matrix $\mathscr{S}_p$ with respect to its enclosing disk. The resulting matrices allow the disk-based multiple scattering formalism to be applied to noncircular particles.

	Specifically, consider an elastic particle $\Omega$ embedded in a homogeneous medium. Given an incident wave $\uinc$, the exterior scattered field  can be represented as:
	\begin{align}\label{intrpre}
	\begin{gathered}
		\uscd(\bx) = \int_{\partial \Omega}[T_{\nu(y)}\Phi(\bx,\by)]^\top\boldsymbol{\varphi}_D(\by)ds_{\by}, \quad \mbox{for the particle with rigid boundary,}\\
		\uscn(\bx) = \int_{\partial \Omega}\Phi(\bx,\by)\boldsymbol{\varphi}_N(\by)ds_{\by}, \quad \mbox{for the particle with traction-free boundary},
		\end{gathered}
	\end{align} 
	where $\bx\in \mr\setminus\overline{\Omega}$. This representation leads to the following second-kind boundary integral equation for the unknown density $\boldsymbol{\varphi}_D$ or $\boldsymbol{\varphi}_N$:
	\begin{align}\label{bdint}
	\begin{gathered}
			\frac{1}{2}\boldsymbol{\varphi}_D+\int_{\partial \Omega}[T_{\nu(y)}\Phi(\bx,\by)]^\top\boldsymbol{\varphi}_D(\by)ds_{\by} = -\uinc,\quad x\in\partial \Omega, \\
			-\frac{1}{2}\boldsymbol{\varphi}_N+\int_{\partial \Omega}[T_{\nu(x)}\Phi(\bx,\by)]^\top\boldsymbol{\varphi}_N(\by)ds_{\by} = -T_{\nu}\uinc,\quad x\in\partial \Omega.
			\end{gathered}
	\end{align}
	\begin{remark}
		We assume that the boundary integral equation \eqref{bdint} is uniquely solvable. If this assumption fails for a particular particle $D_j$, one may instead use a combined-potential formulation \cite{hsiaoBoundaryIntegralEquations2008}.
	\end{remark}
	Equation \eqref{bdint} is discretized by a Nystr{\"o}m method \cite{bao2024singularity,dongHighlyAccurateBoundary2021} based on the piecewise polynomial interpolation scheme of \cite{yaoRobustHighPrecision2024}, which converges exponentially for smooth geometries and remains efficient for domains with corners.
	
	
	In the case of $M$ particles, suppose that particle $D_j$, $1\le j\le M$, is enclosed by a circle $S_j$ centered at $O_j$ with radius $R$. For simplicity, we describe the construction using the single-layer representation of the scattered field. To construct the scattering matrix for $D_j$, we sequentially choose the single-mode basis functions
	\[\nabla(J_n^{(1)}(\kp r)\me^{\mi n \theta}), \quad \nabla^\bot(J_n^{(1)}(\ks r)\me^{\mi n \theta}),\quad n = -p,\cdots,p,\]
	with respect to $O_j$ as incident fields for $D_j$, and solve for the density $\boldsymbol{\varphi}(\by) = (\varphi_1,\varphi_2)^\top$ from the boundary integral equation \eqref{bdint}. We then precompute the multipole expansion coefficients $\{\alpha_l,\beta_l\}$ in \eqref{scatspan} from these source distributions. They are given by
	\begin{align*}
		\begin{aligned}
			\begin{bmatrix}
				\alpha_l\\
				\beta_l
			\end{bmatrix}
			=
			\frac{\mi}{8\omega^2}\int_{\partial D_j}
			\begin{bmatrix}
				\kp \Big(E_{l-1}^{\kp}-E_{l+1}^{\kp}\Big)&
				-\mi \kp\Big(E_{l-1}^{\kp}+E_{l+1}^{\kp}\Big)\\
				\mi\ks \Big(E_{l-1}^{\ks}+E_{l+1}^{\ks}\Big)&
				\ks\Big(E_{l-1}^{\ks}-E_{l+1}^{\ks}\Big)
			\end{bmatrix}
			\begin{bmatrix}
				\varphi_1\\ 	\varphi_2
			\end{bmatrix}	ds_{\by}, \quad
			l = -p,\cdots,p,
		\end{aligned}
	\end{align*}
	where $E_l^k(\by) = J_l(k|\by|)\me^{-\mi l \theta_j(\by)}$. Here, $\by$ denotes the position of a point on $\partial D_j$ relative to $O_j$, and $\theta_j(\by)$ is its polar angle with respect to $O_j$. The formulas for $\alpha_l$ and $\beta_l$ follow from Graf's addition theorem and recurrence relations for Hankel functions and their derivatives \cite{olverNISTHandbookMathematical2010}.
	
	Note that the integral equation \eqref{bdint} needs to be factorized only once, for example by LU factorization or another direct solver, and the factorization can then be applied to the different right-hand sides generated by the incident basis functions. After this computation has been performed for all $-p\leq n\leq p$, we obtain the scattering matrix $\mathscr{S}^j$ for $D_j$. 
	Allowing a different scattering matrix for each particle accounts for the possibility that the particles $D_j$ have distinct geometries. For simplicity, we assume here that the particles have the same shape, or belong to a small number of prototype shapes, up to rotation. In that case, the scattering matrix needs to be computed only once for each prototype particle in its enclosing disk. Matrices for rotated copies are obtained by a simple rotation of the basis. Error analysis for the single-obstacle scattering matrix method in the acoustic setting can be found in \cite{ganeshConvergenceAnalysisParameter2012}.
	
	\begin{remark}
		Several other methods can also be used to evaluate the scattering matrix $\mathscr{S}^j$. One popular approach is the extended boundary condition method (EBCM) \cite{martinMultipleScatteringInteraction2006}, although it may suffer from numerical instabilities for particles with large aspect ratios. Another option is to compute the outgoing coefficients $\{\alpha_l,\beta_l\}$ by the direct projection formula
		\begin{align*}
			\begin{bmatrix}
				\alpha_l\\ \beta_l
			\end{bmatrix} = \frac{1}{2\pi R}\begin{bmatrix}
				\kp H_l^{(1)'}(\kp R) & -\frac{\mi l}{R} H_l^{(1)}(\ks R)\\
				\frac{\mi l}{R}H_l^{(1)}(\kp R)& \ks H_l^{(1)'}(\ks R)
			\end{bmatrix}^{-1}\begin{bmatrix}
				\int_{S_j}\usc\cdot \me^{\mi l\theta}\hat{r}ds\\
				\int_{S_j}\usc\cdot \me^{\mi l\theta}\hat{\theta}ds
			\end{bmatrix},
		\end{align*}
		where $\usc$ is the scattered field on $S_j$ computed from \eqref{intrpre}. This projection formula can also be used to construct the scattering matrix numerically and is commonly adopted in three dimensions \cite{laiFastInverseElastic2022}.
	\end{remark}
	
	Once the scattering matrix is known, the full elastic scattering problem for many particles can be reduced to a multiple scattering problem posed on the enclosing disks. Thus, the multiple scattering method presented in Section \ref{scadisk} applies after replacing the disk scattering matrix by the scattering matrix of each enclosed particle.
	
	This reduction has two main advantages. First, the number of degrees of freedom is substantially reduced because each particle is represented by multipole expansion coefficients rather than by all boundary discretization unknowns. Second, the scattering matrix for each particle is precomputed in isolation, leading to a better-conditioned global multiple scattering system. The resulting system can be solved by combining GMRES with the FMM, reducing the per-iteration cost from $\oo(M^2p^2)$ to $\oo(Mp^2)$.
	
	The main limitations are that a modest separation distance between particles is required and that the bookkeeping becomes more involved when many distinct prototype shapes are present. In many experimental and engineering settings, however, both requirements are naturally satisfied.
	\begin{remark}
		The method above can be viewed as a reduced-order model for the scattering problem. Similar ideas have been successfully applied to acoustic \cite{laiFastSolverMultiparticle2014}, electromagnetic \cite{gimbutasFastMultiparticleScattering2013}, and elastic \cite{laiFastInverseElastic2022,laiFrameworkSimulationMultiple2019} scattering by multiple particles. A complete error analysis is highly involved because of the multiple scattering interactions. Numerically, spectral accuracy is observed as the truncation number $p$ increases. A detailed error analysis for the full multiple scattering problem is beyond the scope of this paper.
	\end{remark}
	\section{Multi-particle scattering in a layered medium}\label{multi_part_layer}
	We have so far discussed the layered background and the multiple scattering formalism separately. We now combine these two components to treat the full problem, in which multiple well-separated particles are embedded in the lower half-space $\lowsp$. By coupling the Sommerfeld representation for the layered medium with the multipole representation for the scattered fields by particles, we write
	\begin{align*}
		\bu_1(\bx) &= \Phi(\bx,\bxs)\bal+\bu_1^s,\\
		\bu_2(\bx) &= \bu_2^s + \sum_{m=1}^{M}\sum_{n=-p}^{p}\left(\alpha_n^m\nabla(H_n^{(1)}(\kp' r_m)\me^{\mi n \theta_m})+\beta_n^m\nabla^\bot(H_n^{(1)}(\ks' r_m)\me^{\mi n \theta_m})\right).
	\end{align*}
	It remains to describe the discretization of the Sommerfeld integral and the construction of the global linear system for the unknown spectral densities $\sg_1,\sg_2$ and multipole coefficients $\{\vec{\alpha}^m,\vec{\beta}^m,\,m=1,2,\cdots,M\}$.
	\subsection{Evaluation of the Sommerfeld integral}
	The numerical evaluation of the Sommerfeld integrals \eqref{us1} and \eqref{us2} has been extensively studied because of their central role in layered medium scattering. We do not attempt to review all available schemes here. Most of them rely on contour deformation in the complex $\xi$-plane to avoid the square-root singularities in the denominator.
	For simplicity, we use the hyperbolic tangent contour \cite{oneilEfficientRepresentationHalfspace2014}
	\[\xi(t) = t-\frac{\mi\tanh(t)}{2}, -t_{max}\leq t\leq t_{max},\]
	for some $t_{max}$.
	If the integrand has decayed to high precision at the endpoints $\pm t_{max}$, which is reasonable when the point source and buried obstacles are not close to the interface $\face$, then the trapezoidal rule yields a spectrally accurate quadrature scheme.
	Let $N_S$ denote the number of quadrature points on the Sommerfeld contour. Each discretization point $\xi_j$ corresponds to a plane-wave mode. We use the same contour and the same set of quadrature points $\{\xi_j\}$ for both $\bu^s_1$ and $\bu^s_2$.
	\subsection{The full linear system}
	Let $\sg$ denote the discretized densities on the interface $\face$, with $\sg = (\sigma_{11},\sigma_{12},\sigma_{21},\sigma_{22})^\top$, and let $\boldsymbol{\gamma}=\{\vec{\alpha}^m,\vec{\beta}^m\}$ collect the multipole coefficients for all $M$ particles in the lower half-space. The vectors $\sigma_{11},\sigma_{12},\sigma_{21},\sigma_{22}$ each have length $N_S$. The full linear system governing multiple scattering in the layered medium can then be written as the block-structured $2\times 2$ system
	\begin{align}\label{aequ}
		\begin{pmatrix}
			A&B\\C&D
		\end{pmatrix}\begin{pmatrix}
			\sg\\ \boldsymbol{\gamma}
		\end{pmatrix}=\begin{pmatrix}
			\mathbf{b}\\ \mathbf{0}
		\end{pmatrix},
	\end{align}	
	where $A$ is a block-diagonal $4N_S\times4N_S$ matrix with $4\times4$ blocks of the form $A_\xi$ in \eqref{aa}, each associated with a quadrature point $\xi_j$ on the contour. The right-hand side component $\mathbf{b}$ is obtained by evaluating the right-hand side of \eqref{ae} at each $\xi_j$. The matrix $D$, defined in \eqref{mul2},
	is the multiple scattering system for the particles.
	The off-diagonal blocks $B$ and $C$ represent the coupling between the particles and the layered interface. More precisely, $B$ maps the particle multipole expansion coefficients to the Sommerfeld representation on the interface $\face$, whereas $C$ evaluates the Sommerfeld integral contributions as incoming local expansions on the scattering disks. We now describe the explicit construction of these two coupling operators.
	\begin{remark}
		We solve \eqref{aequ} using GMRES \cite{saadGMRESGeneralizedMinimal1986} with FMM acceleration. However, the unknowns $\sg$ and $\boldsymbol{\gamma}$ may be poorly scaled relative to each other. Since $A$ is block diagonal with small $4\times4$ blocks, we invert $A$ directly and apply GMRES to the Schur complement of \eqref{aequ}. That is, we solve
	\begin{align*}
		[D-CA^{-1}B][\boldsymbol{\gamma}]=-CA^{-1}\mathbf{b}
	\end{align*}
	instead. This reduced system is better conditioned and involves only the multipole unknowns $\boldsymbol{\gamma}$. The Schur complement also has a natural physical interpretation: it reformulates the scattering problem in terms of the layered Green's function.
	\end{remark}
	\subsection{The Sommerfeld-to-local operator}\label{con1}
	A direct way to map the Sommerfeld density variables $\sg$ to local expansions on each disk is to use the Jacobi--Anger formula \cite{olverNISTHandbookMathematical2010}
	\begin{align}\label{ja}
		\me^{\mi kr\cos\theta}=\sum_{n=-\infty}^\infty \mi^nJ_n(kr)\me^{\mi n\theta}.
	\end{align}
	We first compute the contribution from $\sigma_{21}$ and $\sigma_{22}$ to a local expansion on a disk centered at $(x_1,y_1)$. Using \eqref{ja}, one obtains
	\begin{align}\label{stm}
		\me^{\sqrt{\xi_j^2-{k}^2}y+\mi\xi_j(x-x_0)}=\me^{\sqrt{\xi_j^2-{k}^2}y_1+\mi\xi_j(x_1-x_0)}\sum_{n=-\infty}^\infty \mi^nJ_n(kr)\me^{\mi n(\phi+\theta)},
	\end{align}
	where $\phi = \arccos(\xi_j/k)$, $\theta = \arccos((x-x_1)/r)$, and $r = \sqrt{(x-x_1)^2+(y-y_1)^2}$.
	
	The first column of the spectral Green kernel $\widetilde{\Phi}^2(\xi,\mathbf x,\mathbf{x_0})$ defined in \eqref{lkel} for the lower half-space can be written as
	\begin{align}\label{c1}
		&\frac{1}{4\pi\omega^2}\begin{pmatrix}
			\ks'^2\frac{\me^{\sqrt{\xi_j^2-\ks'^2}y}}{\sqrt{\xi_j^2-\ks'^2}}\me^{\mi\xi_j(x-x_0)}-\xi_j^2\frac{\me^{\sqrt{\xi_j^2-\ks'^2}y}}{\sqrt{\xi_j^2-\ks'^2}}\me^{\mi\xi_j(x-x_0)}+\xi_j^2\frac{\me^{\sqrt{\xi_j^2-\kp'^2}y}}{\sqrt{\xi_j^2-\kp'^2}}\me^{\mi\xi_j(x-x_0)}\\
			\mi\xi_j\me^{\sqrt{\xi_j^2-\ks'^2}y}\me^{\mi\xi_j(x-x_0)}-\mi\xi_j\me^{\sqrt{\xi_j^2-\kp'^2}y}\me^{\mi\xi_j(x-x_0)}
		\end{pmatrix}\nonumber\\
		&=-\frac{\mi\xi_j}{4\pi\omega^2\sqrt{\xi_j^2-\kp'^2}}\nabla\left(\me^{\mi\xi_j(x-x_0)+\sqrt{\xi_j^2-\kp'^2}y}\right)+\frac{1}{4\pi\omega^2}\nabla^\bot\left(\me^{\mi\xi_j(x-x_0)+\sqrt{\xi_j^2-\ks'^2}y}\right).
	\end{align}
	Similarly, its second column can be simplified as
	\begin{align}\label{c2}
		&\frac{1}{4\pi\omega^2}\begin{pmatrix}
			\mi\xi_j\me^{\sqrt{\xi_j^2-\ks'^2}y}\me^{\mi\xi_j(x-x_0)}-\mi\xi_j\me^{\sqrt{\xi_j^2-\kp'^2}y}\me^{\mi\xi_j(x-x_0)}\\
			\ks'^2\frac{\me^{\sqrt{\xi_j^2-\ks'^2}y}}{\sqrt{\xi_j^2-\ks'^2}}\me^{\mi\xi_j(x-x_0)}+{\sqrt{\xi_j^2-\ks'^2}}{\me^{\sqrt{\xi_j^2-\ks'^2}y}}\me^{\mi\xi_j(x-x_0)}-{\sqrt{\xi_j^2-\kp'^2}}\me^{\sqrt{\xi_j^2-\kp'^2}y}\me^{\mi\xi_j(x-x_0)}
		\end{pmatrix}\nonumber\\
		&=-\frac{1}{4\pi\omega^2}\nabla\left(\me^{\mi\xi_j(x-x_0)+\sqrt{\xi_j^2-\kp'^2}y}\right)-\frac{\mi\xi_j}{4\pi\omega^2\sqrt{\xi_j^2-\ks'^2}}\nabla^\bot\left(\me^{\mi\xi_j(x-x_0)+\sqrt{\xi_j^2-\ks'^2}y}\right),
	\end{align}
	These identities imply that $\sigma_{21}$ and $\sigma_{22}$ can be converted to a local expansion on a disk centered at $(x_1,y_1)$ using \eqref{stm}.
	
	Combining \eqref{stm} with \eqref{c1}--\eqref{c2}, the Sommerfeld-to-local conversion formulas are
	\begin{align}\label{eee1}
    \begin{aligned}
        \sigma_{21}(\xi_j)\rightarrow\frac{\sigma_{21}(\xi_j)}{4\pi\omega^2}\me^{\mi\xi_j(x_m-x_0)}\mi^n\begin{pmatrix}
			-\frac{\mi\xi_j}{\sqrt{\xi_j^2-\kp'^2}}\me^{\sqrt{\xi_j^2-{\kp'}^2}y_m+\mi n\phi_{\kp'}},&
			\me^{\sqrt{\xi_j^2-{\ks'}^2}y_m+\mi n\phi_{\ks'}}
		\end{pmatrix}, \\
		\sigma_{22}(\xi_j)\rightarrow\frac{\sigma_{22}(\xi_j)}{4\pi\omega^2}\me^{\mi\xi_j(x_m-x_0)}\mi^n\begin{pmatrix}
			-\me^{\sqrt{\xi_j^2-{\kp'}^2}y_m+\mi n\phi_{\kp'}},&
			-\frac{\mi\xi_j}{\sqrt{\xi_j^2-\ks'^2}}\me^{\sqrt{\xi_j^2-{\ks'}^2}y_m+\mi n\phi_{\ks'}}
		\end{pmatrix},\\
		m = 1,2,\cdots,M,\qquad n = -p,\cdots,p.
    \end{aligned}	
	\end{align}
	Using \eqref{eee1} to apply the off-diagonal block $C$ in \eqref{aequ} costs $\oo(MN_S(4p + 2))$, where $M$ is the number of particles, $N_S$ is the number of Sommerfeld quadrature points, and $p$ is the expansion order in the multiple scattering representation. This cost is acceptable when either $N_S$ or $M$ is modest. For high-frequency problems with many inclusions, where $\kp'$ and $\ks'$ are large and $N_S = \oo(\max\{\kp',\ks'\})$, a more efficient scheme based on the nonuniform FFT (NUFFT) can be used \cite{laiFastSolverMultiparticle2014}. We do not describe that acceleration in detail here.
	\subsection{The multipole-to-Sommerfeld operator}\label{con2}
	The off-diagonal block $B$ in \eqref{aequ} requires a formula that recasts the multipole expansion into the corresponding Sommerfeld representation on the interface $\face$. 
	\begin{lemma}\cite{chengAdaptiveFastSolver2006}
		Let $(x_m,y_m)$ denote the center of a multipole expansion in the lower half-space, with $y_m<0$, and let $(r,\theta)$ denote the polar coordinates of a target point with respect to that center. Then, for $y>y_m$,
		\begin{align}\label{pts1}
			H_n(kr)\me^{\mi n\theta}=\frac{(-1)^n}{\mi\pi}\int_{-\infty}^\infty\frac{\me^{\sqrt{\xi^2-k^2}(y_m-y)}}{\sqrt{\xi^2-k^2}}\me^{\mi\xi(x-x_m)}\left(\frac{\sqrt{\xi^2-k^2}-\xi}{\mi k}\right)^nd\xi.
		\end{align}
	\end{lemma}
	Differentiating both sides of \eqref{pts1} gives
	\begin{align}\label{sed1}
		\nabla(H_n(\kp' r_m)\me^{\mi n\theta_m})
		&=
		\frac{(-1)^n}{\mi\pi}\int_{-\infty}^\infty\begin{pmatrix}
			\mi\xi\\
			-\lp'
		\end{pmatrix}\frac{\me^{\lp'(y_m-y)}}{\lp'}\me^{\mi\xi(x-x_m)}\left(\frac{\lp'-\xi}{\mi\kp'}\right)^nd\xi,
	\end{align}
	\begin{align}\label{sed2}
		\nabla^\bot(H_n(\ks' r_m)\me^{\mi n\theta_m})
		&=\frac{(-1)^n}{\mi\pi}\int_{-\infty}^\infty\begin{pmatrix}
			\ls'\\
			\mi\xi
		\end{pmatrix}\frac{\me^{\ls'(y_m-y)}}{\ls'}\me^{\mi\xi(x-x_m)}\left(\frac{\ls'-\xi}{\mi\ks'}\right)^nd\xi.
	\end{align}
	
	Setting $y=0$ and substituting \eqref{sed1}--\eqref{sed2} into the spectral boundary conditions, a mode-by-mode calculation for each value of $\xi$ gives the contribution of $\vec{\alpha}^m$ and $\vec{\beta}^m$ to the spectral field values:
	\begin{align}\label{b1}
		\begin{aligned}
			{\alpha}_n^m
			\rightarrow\begin{pmatrix}
				{4\omega^2\xi}(-1)^n\frac{\me^{\lp'y_m}}{\lp'}\me^{\mi \xi(x_0-x_m)}\left(\frac{\lp'-\xi}{\mi\kp'}\right)^n\\
				{4\omega^2}{\mi}(-1)^n{\me^{\lp'y_m}}\me^{\mi \xi(x_0-x_m)}\left(\frac{\lp'-\xi}{\mi\kp'}\right)^n
			\end{pmatrix}:=\begin{pmatrix}
				f_{11}(\xi)\\f_{12}(\xi)
			\end{pmatrix},\\
			{\beta}_n^m
			\rightarrow\begin{pmatrix}
				-{4\omega^2}{\mi}(-1)^n{\me^{\ls'y_m}}\me^{\mi \xi(x_0-x_m)}\left(\frac{\ls'-\xi}{\mi\ks'}\right)^n\\
				{4\omega^2\xi}(-1)^n\frac{\me^{\ls'y_m}}{\ls'}\me^{\mi \xi(x_0-x_m)}\left(\frac{\ls'-\xi}{\mi\ks'}\right)^n
			\end{pmatrix}:=\begin{pmatrix}
				f_{21}(\xi)\\f_{22}(\xi)
			\end{pmatrix}.
		\end{aligned}
	\end{align}
	The corresponding contribution to the spectral traction values is
	\begin{align}\label{b2}
		\begin{aligned}
			{\alpha}_n^m
			\rightarrow
			\begin{pmatrix}
				\mu' f_{11}(\xi)(-\lp')+\mu' f_{12}(\xi)(\mi \xi)\\
				(\lambda'+2\mu') f_{12}(\xi)(-\lp')+\lambda' f_{11}(\xi)(\mi \xi)
			\end{pmatrix},~~
			{\beta}_n^m
			\rightarrow
			\begin{pmatrix}
				\mu' f_{21}(\xi)(-\ls')+\mu' f_{22}(\xi)(\mi \xi)\\
				(\lambda'+2\mu') f_{22}(\xi)(-\ls')+\lambda' f_{21}(\xi)(\mi \xi)
			\end{pmatrix}.
		\end{aligned}
	\end{align}
	
	In summary, \eqref{b1} and \eqref{b2} convert the multipole expansion to the corresponding Sommerfeld representation on the interface $\face$. Each multipole coefficient in the expansion of $D_m$ contributes to each of the $N_S$ Sommerfeld quadrature points, so a direct application of the block $B$ costs $\oo((4p + 2)N_SM)$. As for the block $C$, this conversion can also be accelerated using the NUFFT \cite{laiFastSolverMultiparticle2014}.
    

\section{Extension to three-dimensional layered multiple elastic scattering}\label{extension3D}
To demonstrate that the proposed framework extends naturally to three dimensions, we briefly describe its 3D implementation. The main ingredients are the Sommerfeld representation of the 3D elastic Green's tensor, the multipole-to-Sommerfeld and Sommerfeld-to-multipole translation operators, and the numerical discretization of the resulting Sommerfeld integrals. The FMM for 3D elastic scattering in a homogeneous medium has been studied in the work \cite{laiFastInverseElastic2022}, and therefore we do not repeat those details here.

\subsection{Sommerfeld representation of the 3D layered scattering}
The 3D free-space elastic Green's tensor is
                \[
\Phi(\bx, \bxs)=\frac{\kappa_{\mathfrak{s}}^{2}}{4\pi\omega^{2}}\frac{\mathrm{e}^{\mathrm{i}\kappa_{\mathfrak{s}}|\bx-\bxs|}}{|\bx-\bxs|}I+\frac{1}{4\pi\omega^{2}}\nabla\nabla^{\top}\left[\frac{\mathrm{e}^{\mathrm{i}\kappa_{\mathfrak{s}}|\bx-\bxs|}}{|\bx-\bxs|}-\frac{\mathrm{e}^{\mathrm{i}\kappa_{\mathfrak{p}}|\bx-\bxs|}}{|\bx-\bxs|}\right].
\]
Using the Sommerfeld integral representation
$$
\frac{\mathrm{e}^{\mathrm{i}k|\bx-\bxs|}}{4\pi|\bx-\bxs|}
=\iint_{-\infty}^{\infty} \widetilde{g}_{k}(\xi, \mathbf{x}, \mathbf{x}_{0}) d\xi,  \quad \widetilde{g}_{k}(\xi, \mathbf{x}, \mathbf{x}_{0})=\frac{\mi}{8\pi^2}\frac{\me^{\mi k_{z} |z - z_0|}}{k_{z}} \me^{\mi \xi\cdot(r_{xy}-r_{xy}')},
$$
\[\xi=(k_x,k_y), r_{xy}-r_{xy}'=(x-x_0,y-y_0), k_z=\sqrt{k^2-k_x^2-k_y^2},\]
we obtain the corresponding spectral representation
\[
\widetilde{\Phi}\left(\xi, \mathbf{x}, \mathbf{x_{0}}\right)=\frac{1}{\omega^{2}}\left(\kappa_{\mathfrak{s}}^{2} \widetilde{g}_{\kappa_{\mathfrak{s}}} \mathbb{I}+\nabla \nabla^{\top}\left(\widetilde{g}_{\kappa_{\mathfrak{s}}}-\widetilde{g}_{\kappa_{\mathfrak{p}}}\right)\right),\quad 
\Phi(\mathbf{x}, \mathbf{x}_0) = \iint_{-\infty}^{\infty} \widetilde{\Phi}(\xi, \mathbf{x}, \mathbf{x}_0) d\xi.
\]
\begin{remark}
	Numerically, the 3D Sommerfeld integral can be evaluated in polar coordinates $(k_\rho,k_\theta)$, with a contour deformation applied to the radial variable \cite{cuiFastEvaluationSommerfeld1999}:
\begin{align}\label{condd1}
    k_\rho=\lambda,\lambda(t) = t - \mi\alpha\tanh(\beta t), \quad t\in [0,T_{max}].
\end{align}
In particular, we use $N_r$ Gauss--Legendre points in the radial direction and $N_\theta$ equally spaced points in the angular direction. It holds   
\begin{align*}
   \frac{\mathrm{e}^{\mathrm{i}k|\bx-\bxs|}}{4\pi|\bx-\bxs|}
    &\approx\frac{\mi}{8\pi^2}\sum_{i=1}^{N_r}  \frac{\omega_i\lambda(t_i)\lambda'(t_i)}{\sqrt{k^2-\lambda_i^2}} \me^{\mi \sqrt{k^2-\lambda_i^2} |z - z_0|}\frac{2\pi}{N_\theta}\sum_{j=1}^{N_\theta}\me^{\mi\lambda(t_i)\cos (\theta_j-\phi)}\\
    &=\frac{\mi}{4\pi N_\theta}\sum_{i=1}^{N_r}  \frac{\omega_i\lambda(t_i)\lambda'(t_i)}{\sqrt{k^2-\lambda_i^2}}  \sum_{j=1}^{N_\theta}\me^{\mi\lambda(t_i)\cos\theta_j(x-x_0)+\mi\lambda(t_i)\sin\theta_j(y-y_0)+\mi{\sqrt{k^2-\lambda_i^2}} |z - z_0|}.
\end{align*}
\end{remark}

The single spectral layer potentials in the upper and lower half-spaces are given by
\[
\usc=\iint_{-\infty}^{\infty}\widetilde{\Phi}^{1}\left(\xi,\mathbf{x},\mathbf{x}_{0}\right)\sg(\xi) d\xi,\qquad \uscp=\iint_{-\infty}^{\infty}\widetilde{\Phi}^{2}\left(\xi,\mathbf{x},\mathbf{x}_{0}\right)\sg'(\xi) d\xi,
\]
where 
\begin{align*}
    \widetilde{\Phi}^{1}(\xi, \mathbf{x}, \mathbf{x_{0}})=\frac{1}{\omega^{2}}\left(\kappa_{\mathfrak{s}}^{2} \widetilde{g}_{\kappa_{\mathfrak{s}}}^{1} \mathbb{I}+\nabla \nabla^{\top}(\widetilde{g}_{\kappa_{\mathfrak{s}}}^{1}-\widetilde{g}_{\kappa_{\mathfrak{p}}}^{1})\right), \quad \widetilde{g}_{k}^{1}(\xi, \mathbf{x}, \mathbf{x_{0}})=\frac{\mi}{8\pi^2}\frac{\me^{\mi k_{z} z}}{k_{z}} \me^{\mi  \xi \cdot(r_{xy}-r_{xy}')},\\
    \widetilde{\Phi}^{2}(\xi, \mathbf{x}, \mathbf{x_{0}})=\frac{1}{\omega^{2}}\left(\kappa_{\mathfrak{s}}'^{2} \widetilde{g}_{\kappa_{\mathfrak{s}}}^{2} \mathbb{I}+\nabla \nabla^{\top}(\widetilde{g}_{\kappa_{\mathfrak{s}}}^{2}-\widetilde{g}_{\kappa_{\mathfrak{p}}}^{2})\right), \quad \widetilde{g}_{k}^{2}(\xi, \mathbf{x}, \mathbf{x_{0}})=\frac{\mi}{8\pi^2}\frac{\me^{-\mi k_z' z}}{k_z'} \me^{\mi  \xi \cdot(r_{xy}-r_{xy}')}.
\end{align*}
For point-source incidence, the total field in the upper half-space $(0<z<z')$, consisting of the incident point-source field and the upward layered potential, is
\begin{align*}
    \boldsymbol{(u_1^t)}_{\fp}
    &=\frac{\mi}{8\pi^2\omega^2}\iint_{-\infty}^{\infty} -\nabla\nabla^\top \left[\frac{\me^{\mi \kpz  (z'-z) }}{\kpz } \me^{\mi  \xi \cdot(r_{xy}-r_{xy}')}\right]\boldsymbol{p} - \nabla\nabla^\top \left[\frac{\me^{\mi \kpz  z }}{\kpz } \me^{\mi  \xi \cdot(r_{xy}-r_{xy}')} \right]\boldsymbol{\sigma}(k_x,k_y) dk_x dk_y\\
    &=\frac{\mi\kp^2}{8\pi^2\omega^2}\iint_{-\infty}^{\infty} \me^{\mi  \xi \cdot(r_{xy}-r_{xy}')} \left[\frac{\me^{\mi \kpz  (z'-z) }}{\kpz } \boldsymbol{e_{\fp,1}}\boldsymbol{e_{\fp,1}}^{\top}\boldsymbol{p} +\frac{\me^{\mi \kpz  z }}{\kpz }  \boldsymbol{e_{\fp,2}}\boldsymbol{e_{\fp,2}}^{\top}\boldsymbol{\sigma}(k_x,k_y) \right]dk_x dk_y,
\end{align*}
where $\boldsymbol{e_{\fp,1}}=(k_x,k_y,-\sqrt{\kp^2-k_x^2-k_y^2})^\top/\kp, \boldsymbol{e_{\fp,2}}=(k_x,k_y,\sqrt{\kp^2-k_x^2-k_y^2})^\top/\kp,$ and
\begin{align*}
    \boldsymbol{(u_1^t)}_{\fs}
    &=\frac{\mi\ks^2}{8\pi^2\omega^2}\iint_{-\infty}^{\infty}  \me^{\mi  \xi \cdot(r_{xy}-r_{xy}')}\left[\frac{\me^{\mi \ksz  (z'-z) }}{\ksz } (\mathbb{I}-\boldsymbol{e_{\fs,1}}\boldsymbol{e_{\fs,1}}^{\top})\boldsymbol{p} + \frac{\me^{\mi \ksz  z }}{\ksz }  (\mathbb{I}-\boldsymbol{e_{\fs,2}}\boldsymbol{e_{\fs,2}}^{\top})\boldsymbol{\sigma}(k_x,k_y)\right] dk_x dk_y,
\end{align*}
where $\boldsymbol{e_{\fs,1}}=(k_x,k_y,-\sqrt{\ks^2-k_x^2-k_y^2})^\top/\ks, \boldsymbol{e_{\fs,2}}=(k_x,k_y,\sqrt{\ks^2-k_x^2-k_y^2})^\top/\ks.$

Similarly, the total field in the lower half-space $(z<0)$, consisting of the scattered field from the buried obstacles and the downward layered potential, is
\begin{align*}
   \boldsymbol{u_2^t} =&  (\uscp)_{\fp}+(\uscp)_{\fs}
+ \boldsymbol{u^{s}_{par}}\\
    =&\frac{\mi\kp'^2}{8\pi^2\omega^2}\iint_{-\infty}^{\infty} \me^{\mi  \xi \cdot(r_{xy}-r_{xy}')} \frac{\me^{-\mi \kpz' z }}{\kpz'}  \boldsymbol{e_{\fp',1}}\boldsymbol{e_{\fp',1}}^{\top}\sg'(k_x,k_y) dk_x dk_y\\
&+\frac{\mi\ks'^2}{8\pi^2\omega^2}\iint_{-\infty}^{\infty}  \me^{\mi  \xi \cdot(r_{xy}-r_{xy}')}\frac{\me^{-\mi \ksz' z }}{\ksz'}  (\mathbb{I}-\boldsymbol{e_{\fs',1}}\boldsymbol{e_{\fs',1}}^{\top})\boldsymbol{\sigma}'(k_x,k_y) dk_x dk_y\\
&+\sum_{n=1}^{\infty} \sum_{m=-n}^{n} \left( \alpha_{n, m} \nabla \times \nabla \times \left( x v_{n, m}^{\kappa_{\mathfrak{s}}'}\right) /\left( \mi \kappa_{\mathfrak{s}}'\right) + \beta_{n, m} \nabla \times\left( x v_{n, m}^{\kappa_{\mathfrak{s}}'}\right) \right)+\sum_{n=0}^{\infty} \sum_{m=-n}^{n} \gamma_{n, m} \nabla v_{n, m}^{\kappa_{\mathfrak{p}}'},
\end{align*}
where $\boldsymbol{e_{\fp',1}}=(k_x,k_y,-\sqrt{\kp'^2-k_x^2-k_y^2})^\top/\kp'$ and $\boldsymbol{e_{\fs',1}}=(k_x,k_y,-\sqrt{\ks'^2-k_x^2-k_y^2})^\top/\ks'$. The functions $v_{n, m}^\kappa(r,\theta,\phi)=h_n^{(1)}(\kappa r)Y_n^m(\theta,\phi)$ are spherical wave functions. 

In the 3D elastic case, the traction operator is
\[
T_{\nu}\boldsymbol{u} = 2\mu(\nu\cdot\nabla)\boldsymbol{u} + \lambda\nu(\nabla\cdot\boldsymbol{u}) + \mu\nu\times(\text{curl}\,\boldsymbol{u}).
\]
When $\nu=(0,0,1)^\top$ on the interface, this operator becomes
\begin{align*}
    T_{\nu}\boldsymbol{u} = \begin{pmatrix} 
\mu (\partial_z u_x + \partial_x u_z) \\ 
\mu (\partial_z u_y + \partial_y u_z) \\ 
\lambda (\partial_x u_x + \partial_y u_y) + (\lambda + 2\mu) \partial_z u_z 
\end{pmatrix}.
\end{align*}
Thus, the tractions of $\boldsymbol{u_1^t}$ and $\boldsymbol{u_2^t}$ can also be evaluated in Sommerfeld integral form. Enforcing continuity of displacement and traction across the interface, together with the Dirichlet or Neumann boundary condition on the particles, yields a block-structured $2\times2$ linear system analogous to \eqref{aequ}:
	\begin{align}\label{aequ2}
		\begin{pmatrix}
			A&B\\C&D
		\end{pmatrix}\begin{pmatrix}
			\sg_t\\ \boldsymbol{\tau}
		\end{pmatrix}=\begin{pmatrix}
			\mathbf{b}\\ \mathbf{0}
		\end{pmatrix},
	\end{align}	
where $\sg_t=(\sg,\sg')$ denotes the unknown spectral density on the interface, and $\boldsymbol{\tau}=(\boldsymbol{\alpha},\boldsymbol{\beta},\boldsymbol{\gamma})$ collects the unknown outgoing coefficients on the particles. As in the 2D case, $A$ and $D$ represent the interface-to-interface and particle-to-particle self-interactions, respectively. The matrix $A$ is constructed from the explicit Sommerfeld expressions above and the mode-matched continuity conditions at the interface. The matrix $D$ follows from the 3D homogeneous medium multiple scattering formulation in \cite{laiFastInverseElastic2022}. It remains to describe the construction of the coupling blocks $B$ and $C$.
\subsection{Translation from spherical wave functions to plane-wave factors}
In this subsection, we translate the three spherical-wave components, associated with $\frac{1}{\mi \ks'}\nabla\times\nabla\times$, $\nabla\times$, and $\nabla$, into Sommerfeld representations. A direct vector calculation gives the following representations for the three types of spherical basis functions:
\begin{align*}
\nabla\times\nabla\times (\boldsymbol{r}v_{n,m}^{\ks'}(r,\theta,\phi))/(\mi \ks')
&=\frac{(-\mi)^{n}}{ 2\pi \ks'} \iint_{-\infty}^{\infty}\frac{\me^{\mi\ks'\boldsymbol{e_{\fs',2}}\cdot(\mathbf{x-x}_l)}}{\ksz'} \operatorname{Grad}\mathcal{Y}_n^m(\boldsymbol{e_{\fs',2}})dk_x dk_y,\\
\nabla\times (\boldsymbol{r}v_{n,m}^{\ks'}(r,\theta,\phi))
&=\frac{(-\mi)^{n}}{ 2\pi \ks'} \iint_{-\infty}^{\infty}\frac{\me^{\mi\ks'\boldsymbol{e_{\fs',2}}\cdot(\mathbf{x-x}_l)}}{\ksz'} \operatorname{Grad}\mathcal{Y}_n^m(\boldsymbol{e_{\fs',2}})\times \boldsymbol{e_{\fs',2}}dk_x dk_y,\\
\nabla v_{n,m}^{\kp'}
        &=\frac{(-\mi)^{n-1}}{2\pi} \iint_{-\infty}^{\infty}\frac{\me^{\mi\kp' \boldsymbol{e_{\fp',2}}\cdot(\mathbf{x-x}_l)}}{\kpz'}\mathcal{Y}_n^m(\boldsymbol{e_{\fp',2}})\boldsymbol{e_{\fp',2}} dk_x dk_y,
\end{align*}
where $\boldsymbol{e_{\fp',2}}=(k_x,k_y,\sqrt{\kp'^2-k_x^2-k_y^2})^\top/\kp'$ and $\boldsymbol{e_{\fs',2}}=(k_x,k_y,\sqrt{\ks'^2-k_x^2-k_y^2})^\top/\ks'$. Here, $\mathcal{Y}_n^m$ and $\operatorname{Grad}\mathcal{Y}_n^m$ denote the analytic continuations of the spherical harmonic $Y_n^m$ and its surface gradient $\operatorname{Grad}{Y}_n^m$, respectively. Removing the common plane-wave factor $\me^{\mi \xi\cdot (r_{xy}-r_{xy}')}$ from the formulas above gives the entries of the block $B$.

\subsection{Translation from plane-wave factors to spherical wave functions}
In 3D elastic scattering, the incident wave can be expanded as
        \[
\boldsymbol{u}^{i}(x)=\sum_{n=1}^{\infty} \sum_{m=-n}^{n} \left( a_{n, m} \nabla \times \nabla \times \left( x u_{n, m}^{\kappa_{\mathfrak{s}}'}\right) /\left( \mi \kappa_{\mathfrak{s}}'\right) + b_{n, m} \nabla \times\left( x u_{n, m}^{\kappa_{\mathfrak{s}}'}\right) \right)+\sum_{n=0}^{\infty} \sum_{m=-n}^{n} c_{n, m} \nabla u_{n, m}^{\kappa_{\mathfrak{p}}'}.
\]
For the incident plane wave \[
\boldsymbol{u}^{\text {plane }}(\boldsymbol{x})= \mathrm{e}^{\mathrm{i} \kappa_{\mathrm{s}}' \boldsymbol{x} \cdot \boldsymbol{d}}(\boldsymbol{d} \times \boldsymbol{p}) \times \boldsymbol{d}+ \mathrm{e}^{\mathrm{i} \kappa_{\mathfrak{p}}' \boldsymbol{x} \cdot \boldsymbol{d}}(\boldsymbol{d} \cdot \boldsymbol{p}) \boldsymbol{d},
\] the corresponding coefficients are
\begin{align*}
    a_{n, m}&=\frac{4 \pi \mathrm{i}^{n}}{n(n+1)} \operatorname{Grad} Y_{n}^{-m}(\boldsymbol{d}) \cdot \boldsymbol{p},\\
    b_{n, m}&= -\frac{4 \pi \mathrm{i}^{n}}{n(n+1)} \boldsymbol{d} \times \operatorname{Grad} Y_{n}^{-m}(\boldsymbol{d}) \cdot \boldsymbol{p},\\
    c_{n, m}&=-\frac{4 \pi \mathrm{i}^{n+1}}{\kappa_{\mathfrak{p}}} Y_{n}^{-m}(\boldsymbol{d}) \boldsymbol{d} \cdot \boldsymbol{p}.
\end{align*}
The spectral representation can be interpreted as an integral of plane-wave factors over different propagation directions:
\begin{align*}
    (\uscp)_{\fp'}
    &=\frac{\mi}{8\pi^2\omega^2}\iint_{-\infty}^{\infty} \frac{\kp'^2 }{{\kpz'}}\left[\me^{\mi \kp' \boldsymbol{e_{\fp',1}} \cdot(x-x',y-y',z)} (\boldsymbol{e_{\fp',1}}\cdot \boldsymbol{\sigma}(k_x,k_y) )\boldsymbol{e_{\fp',1}}\right] dk_x dk_y,\\
    (\uscp)_{\fs'}
    &=\frac{\mi}{8\pi^2\omega^2}\iint_{-\infty}^{\infty} \frac{\ks'^2}{\ksz'} \left[\me^{\mi \ks'\boldsymbol{e_{\fs',1}} \cdot(x-x',y-y',z)}  \boldsymbol{e_{\fs',1}}\times\boldsymbol{\sigma}(k_x,k_y)\times \boldsymbol{e_{\fs',1}}\right] dk_x dk_y,
\end{align*}
where $\boldsymbol{e_{\fp',1}}=(k_x,k_y,-\sqrt{\kp'^2-k_x^2-k_y^2})^\top/\kp'$ and $\boldsymbol{e_{\fs',1}}=(k_x,k_y,-\sqrt{\ks'^2-k_x^2-k_y^2})^\top/\ks'$. Hence $(\uscp)_{\fp'}$ can be discretized as a sum of compressional plane waves and $(\uscp)_{\fs'}$ as a sum of shear plane waves. Each plane-wave contribution is then translated into the corresponding spherical wave functions, which yields the matrix $C$.

	\section{Numerical experiments}\label{numer_exp}
	In this section, we demonstrate the performance of the proposed algorithm through several numerical examples. For simplicity, the two-dimensional examples use a single family of particle shapes, parameterized by
	\begin{align*}
		x(t) = (a_1+a_2\cos(a_3t))(\cos t,\sin t), \qquad t\in[0,2\pi],
	\end{align*}
	whereas the three-dimensional examples use spherical particles. For fixed $a_1,a_2,a_3$, multiple copies of the particle are randomly distributed and well separated in the lower half-space with random orientations. As discussed in Section \ref{arbq}, particles with more complicated boundaries do not introduce any essential difficulty for the proposed scheme, although the precomputation of the scattering matrix becomes more involved. 

	Throughout the numerical experiments, unless otherwise stated, we set the angular frequency and Lam\'e parameters to $\omega = 2\pi,\lambda = 1,\mu = 2, \lambda' = 3,\mu' = 4$.
	These choices give the compressional wavenumber $\kp=\frac{2\sqrt{5}\pi}{5}$ and shear wavenumber $\ks=\frac{2\sqrt{2}\pi}{2}$ in the upper half-space, and $\kp' = \frac{2\sqrt{11}\pi}{11}$ and $\ks' = \pi$ in the lower half-space. We use the Nystr\"om discretization to construct the scattering matrix and accurately handle the singular integrals \cite{yaoRobustHighPrecision2024}. To choose an appropriate value of $N_S$ for the Sommerfeld discretization, we assume that both the point source defining the incident field and the particle nearest to the interface $\face$ are at least half wavelength away from the interface. Under these assumptions, we set $t_{max}=20$ and discretize
	the integral contour using 200 equispaced points. Since the number of unknowns is large, the linear systems \eqref{aequ} and \eqref{aequ2} are not assembled explicitly. Instead, they are solved by GMRES, with particle-to-particle interactions accelerated by the FMM in both two and three dimensions. All experiments were implemented in MATLAB and carried out on a laptop with an Intel CPU and 16 GB of memory.
	
	The following notation is used in Tables \ref{t1}-\ref{t4}:
	\begin{itemize}
		\item $\omega$: the angular frequency.
		\item $N_{pts}$: the number of discretization points on a single particle.
		\item $N_{particle}$: the total number of particles. 
		\item $N_{term}$: the highest order used in the local and multipole expansions of a single particle, so each expansion has $2N_{term}+1$ terms. 
		\item $N_{tot}$: the total number of particle unknowns. The total number of unknowns in the full linear system also includes the additional $4N_S$ interface unknowns on $\face$. For layer-potential unknowns, $N_{tot}=2N_{pts}N_{particle}$. For multipole expansion unknowns, $N_{tot}=2(2N_{term}+1)N_{particle}$. 
		\item $N_{iter}$: the number of GMRES iterations. 
		\item $T_{solve}$: the GMRES solution time in seconds. 
		\item $E_{err}$: the relative $L^2$ error of the elastic field measured at 12 randomly distributed points in the lower half-space.
	\end{itemize}
	 
	\subsection{Scattering with an analytic solution}
\begin{figure}[!h]
		\centering
		\subfloat[]{
			\includegraphics[scale=0.25]{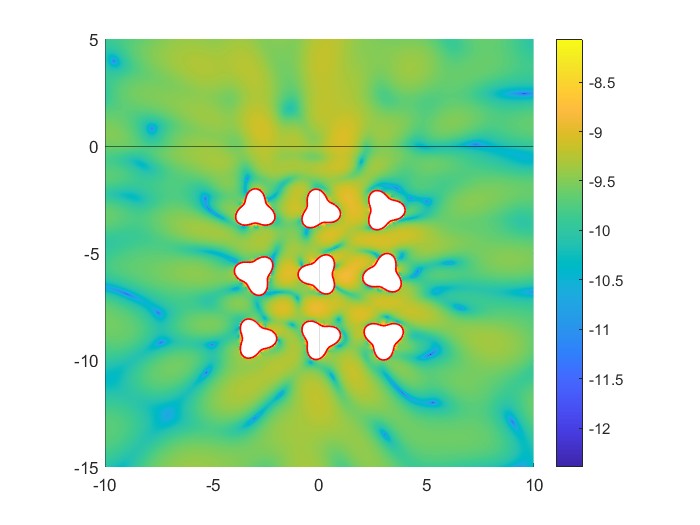}}
		\subfloat[]{
			\includegraphics[scale=0.25]{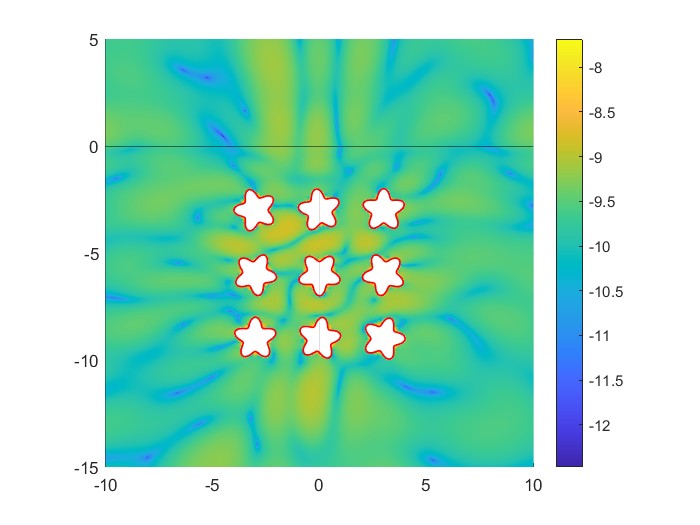}}
		\\
		\subfloat[]{
			\includegraphics[scale=0.23]{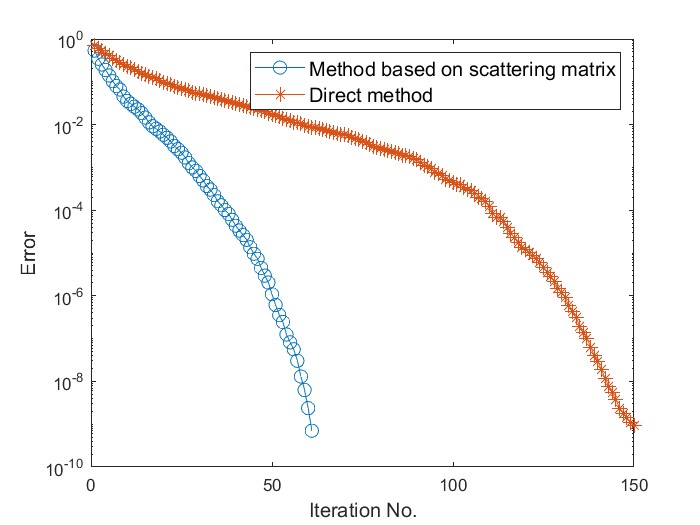}}
		\subfloat[]{
			\includegraphics[scale=0.23]{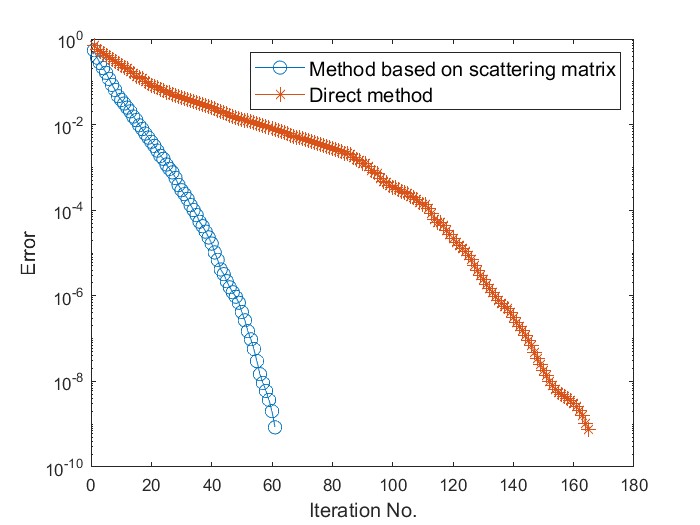}}
		\caption{Elastic scattering of 9 particles at $\omega = 2\pi$. (a) The logarithmic error of the computed field when $a_1 = 4/5, a_2 = 1/5, a_3 = 3$. (b) The logarithmic error of the computed field when $a_1 = 4/5, a_2 = 1/5, a_3 = 5$. (c) Comparison of the GMRES convergence rate when $a_1 = 4/5, a_2 = 1/5, a_3 = 3$. (d) Comparison of the GMRES convergence rate when $a_1 = 4/5, a_2 = 1/5, a_3 = 5$.}
		\label{a1}
	\end{figure}

	In this example, we consider elastic scattering by 9 rigid particles, denoted by $D_i$, $i=1,\cdots,9$. See Figure \ref{a1}(a) and (b) for the illustration. Two methods are compared. The first is a direct method, in which all particle boundaries are discretized by points and the Nystr\"om discretization is applied directly to the integral equation.  In this direct method, we still use the free-space Green's function and treat the interactions between the interface $\face$ and the layer potentials on $\partial D_i$, $i=1,\cdots,9$, as off-diagonal blocks, because computing the layered Green's function directly is expensive. The source-to-Sommerfeld and Sommerfeld-to-source operators are omitted here. They follow the same idea as in Subsections \ref{con1} and \ref{con2}. The second method is the proposed scattering-matrix-based method, in which we first construct the scattering matrix and then solve for the multipole expansion coefficients. To verify the accuracy of both methods, we construct an artificial solution by letting the fields in $\upsp$ and $\lowsp$ be generated by point sources located in $\lowsp$ and inside one buried particle, respectively. Specifically, we choose the exterior elastic field as
	\begin{align*}
		\bu(\bx) = \begin{cases*}
			\nabla H_0^{(1)}(\kp|\bx-\bx_0|),\quad \bx\in\upsp,\\
			\nabla H_0^{(1)}(\kp'|\bx-\bx_0'|),\quad \bx\in\lowsp\setminus\bigcup_{j=1}^9 D_j,
		\end{cases*}
	\end{align*}
	where $\bx_0 = (0,-3)$ and $\bx_0'=(-3,-8.9)$, with $\bx_0'$ located inside $D_7$. By uniqueness, this solution can be recovered by enforcing boundary conditions on $\partial D_j$, $j=1,\cdots,9$, that are consistent with the prescribed field $\bu$.
	
\begin{figure}[!h]
		\centering
		\subfloat[]{
			\includegraphics[scale=0.25]{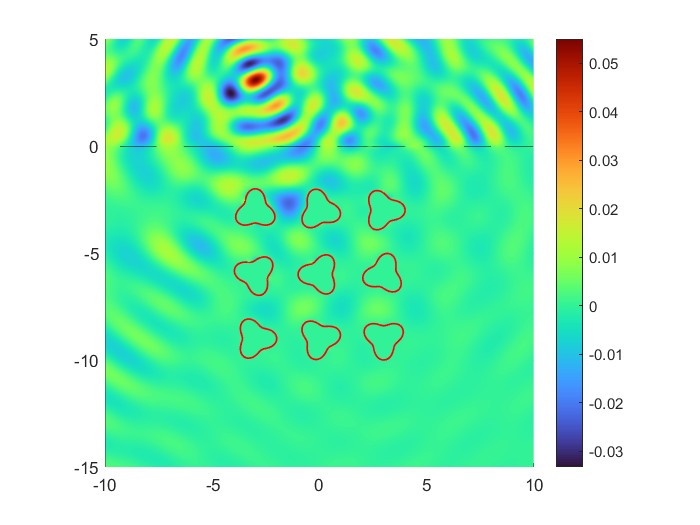}}
		\subfloat[]{
			\includegraphics[scale=0.25]{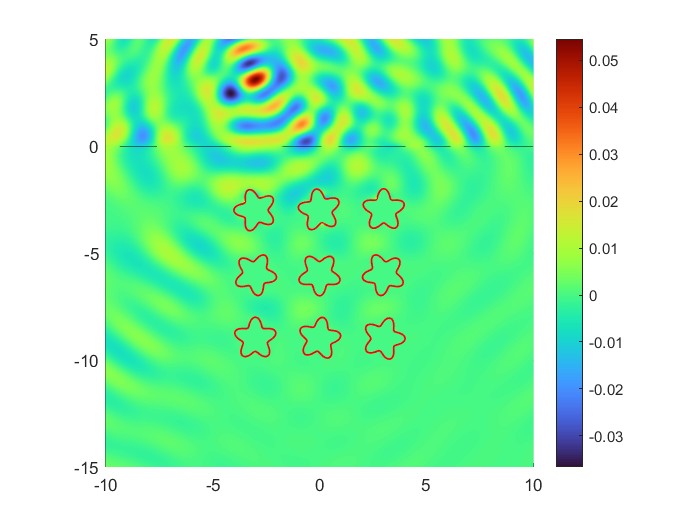}}
		\\
		\subfloat[]{
			\includegraphics[scale=0.25]{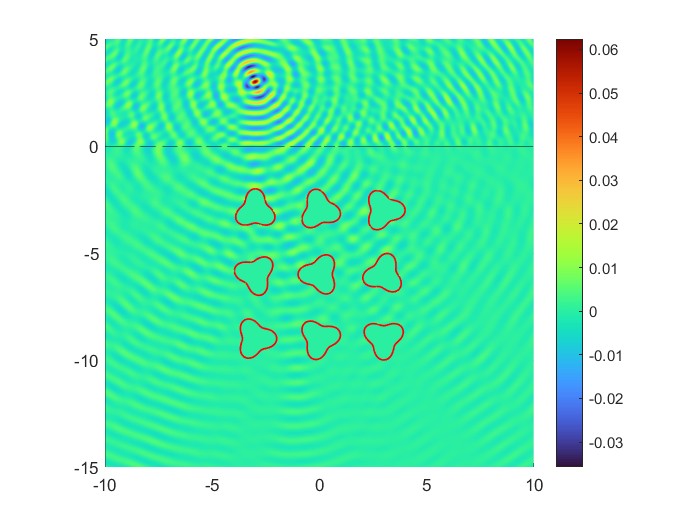}}
		\subfloat[]{
			\includegraphics[scale=0.25]{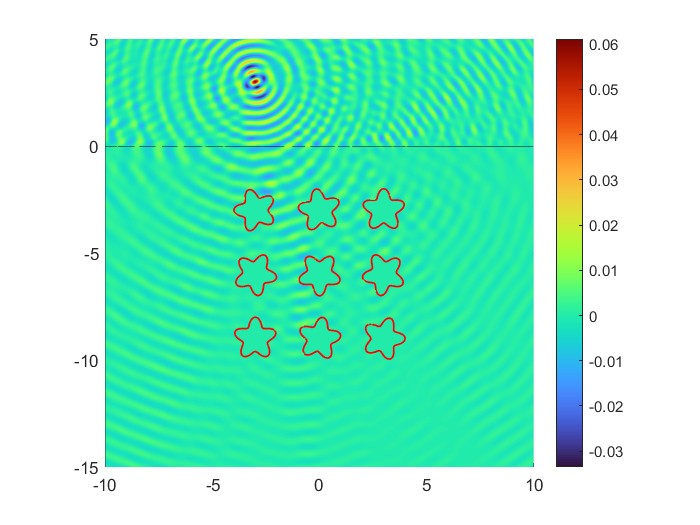}}
		\caption{Elastic scattering of 9 particles by point source illumination. (a) $\omega = 2\pi, a_1 = 4/5, a_2 = 1/5, a_3 = 3$. (b) $\omega = 2\pi, a_1 = 4/5, a_2 = 1/5, a_3 = 5$. (c) $\omega = 6\pi, a_1 = 4/5, a_2 = 1/5, a_3 = 3$. (d) $\omega = 6\pi, a_1 = 4/5, a_2 = 1/5, a_3 = 5$.}
		\label{a9}
	\end{figure}

	We test noncircular geometries with $a_2\neq0$ and $a_3\neq0$. Figure \ref{a1}(a)(b) shows the logarithmic error of the elastic field computed by the scattering matrix based method with $N_{term}=10$ and GMRES tolerance $10^{-9}$, compared with the analytic solution for $\omega=2\pi$.
	The results show more than eight digits of accuracy. Figure \ref{a1}(c)(d) shows that, when each particle is discretized with 64 points in the direct method, nearly 200 GMRES iterations are required for convergence, whereas the scattering matrix based method converges much faster for the same accuracy. Tables \ref{t1} and \ref{t2} present the corresponding results for the direct method and the scattering matrix based method, respectively. With the same number of boundary discretization points in the direct method, the accuracy for starfish-shaped particles is lower than that for pear-shaped particles, whereas the scattering matrix based method is nearly insensitive to the particle shape and achieves high accuracy more quickly. In particular, excluding the precomputation cost of the scattering matrix, the proposed solver is nearly 50 times faster than the direct method for the same accuracy. Figure \ref{a9} shows the imaginary part of the first component of the total elastic field scattered by a point source located at $(-3,3)$ with polarization direction $(\cos(\frac{\pi}{4}),\sin(\frac{\pi}{4}))^\top$ for different particle shapes and frequencies.
	\begin{table}[h!]
		\begin{center}
			\caption{Results for the elastic scattering of 9 particles based on the direct method}
			\label{t1}
			\begin{tabular}{c|cc|ccc|ccc}
				\hline
				& & & $a_1 = \frac{4}{5}$&$a_2=\frac{1}{5}$&$a_3=3$&$a_1 = \frac{4}{5}$&$a_2= \frac{1}{5}$&$a_3=5$\\
				\hline
				$\omega$&$N_{pts}$&$N_{tot}$&$N_{iter}$&$T_{solve}$&$E_{err}$&$N_{iter}$&$T_{solve}$&$E_{err}$\\
				\hline
				\multirow{3}{*}{$\pi$} &16&288&77&2.44E-1&1.86E-2&79&3.53E-1&4.57E-2\\
				&32&576&97&5.39E-1&5.00E-4&89&4.67E-1&1.31E-3\\
				&64&1152&105&9.30E-1&1.80E-7&115&1.01E0&1.55E-5\\
				\hline
				\multirow{3}{*}{$2\pi$} &16&288&128&4.79E-1&4.92E-2&132&4.32E-1&8.03E-2\\
				&32&576&145&7.88E-1&1.97E-3&138&7.26E-1&3.77E-3\\
				&64&1152&150&1.36E0&3.08E-7&165&1.54E0&3.42E-5\\
				\hline
				\multirow{3}{*}{$4\pi$} &16&288&242&9.67E-1&1.51E0&227&9.93E-1&5.54E-1\\
				&32&576&278&1.84E0&1.93E-2&263&1.75E0&3.67E-2\\
				&64&1152&286&3.46E0&1.36E-6&289&3.58E0&1.49E-4\\
				\hline
			\end{tabular}
		\end{center}
	\end{table}
	
	\begin{table}[h!]
		\begin{center}
			\caption{Results for the elastic scattering of 9 particles based on scattering matrix method}
			\label{t2}
			\begin{tabular}{c|cc|ccc|ccc}
				\hline
				& & & $a_1 = \frac{4}{5}$&$a_2=\frac{1}{5}$&$a_3=3$&$a_1 = \frac{4}{5}$&$a_2= \frac{1}{5}$&$a_3=5$\\
				\hline
				$\omega$&$N_{term}$&$N_{tot}$&$N_{iter}$&$T_{solve}$&$E_{err}$&$N_{iter}$&$T_{solve}$&$E_{err}$\\
				\hline
				\multirow{3}{*}{$\pi$} 
				&5&198&43&1.46E-2&1.71E-5&43&1.59E-2&1.52E-5\\
				&10&378&43&4.51E-2&5.32E-9&43&5.64E-2&3.56E-9\\
				&15&558&43&8.24E-2&3.05E-10&43&8.95E-2&3.24E-10\\
				\hline
				\multirow{3}{*}{$2\pi$} 
				&5&198&61&2.21E-2&4.46E-5&61&2.70E-2&4.11E-5\\
				&10&378&61&7.37E-2&8.31E-9&61&8.79E-2&8.38E-9\\
				&15&558&61&1.52E-1&8.84E-10&61&1.32E-1&1.28E-9\\
				\hline
				\multirow{3}{*}{$4\pi$} 
				&5&198&75&2.89E-2&1.23E-3&87&3.93E-2&1.07E-3\\
				&10&378&86&1.56E-1&3.42E-8&93&1.76E-1&4.25E-8\\
				&15&558&86&2.31E-1&1.56E-9&93&2.88E-1&1.92E-9\\
				\hline
			\end{tabular}
		\end{center}
	\end{table}
	
	\subsection{Scattering by point source incidence}
	In this example, we test the algorithm for a large number of buried particles under point source incidence. The point source is located at $(-3,3)$ with polarization direction $(\cos(\frac{\pi}{4}),\sin(\frac{\pi}{4}))^\top$, and all particles are randomly placed in the square domain $[-5,5]\times[-13,-3]$.
	We construct the scattering matrix by solving the integral equation \eqref{bdint} on a single particle with 64 discretization points. The number of terms in the multipole expansion is chosen to be $N_{term} = 5$ 
	and the GMRES tolerance is set to $10^{-6}$. 
	To verify the accuracy of the computed solution, we compare it with the solution obtained using $N_{term} = 20$. Numerical results for scattering by multiple starfish-shaped particles at various angular frequencies are shown in Table \ref{t3}. The table indicates that, for a fixed number of particles, the number of GMRES iterations grows roughly linearly with the angular frequency $\omega$. 
	For high-frequency problems with Neumann boundary conditions, the observed accuracy is lower, indicating that a larger truncation number $N_{term}$ and a tighter GMRES tolerance should be used in such cases.
	
\begin{figure}[!t]
		\centering
		\subfloat[]{
			\includegraphics[scale=0.25]{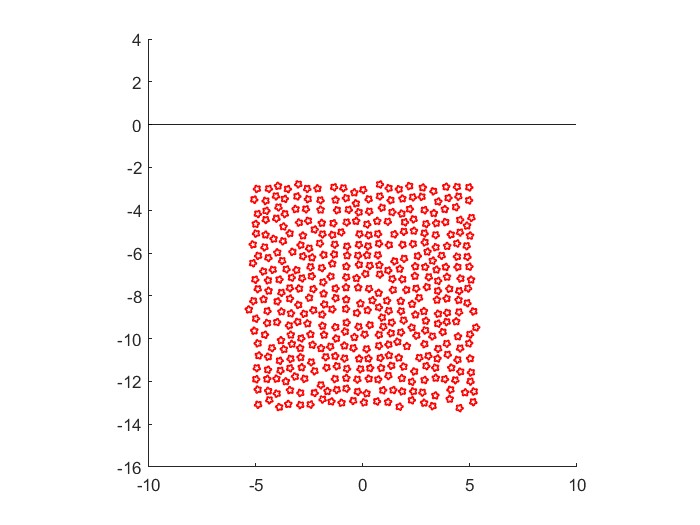}}
		\subfloat[]{
			\includegraphics[scale=0.25]{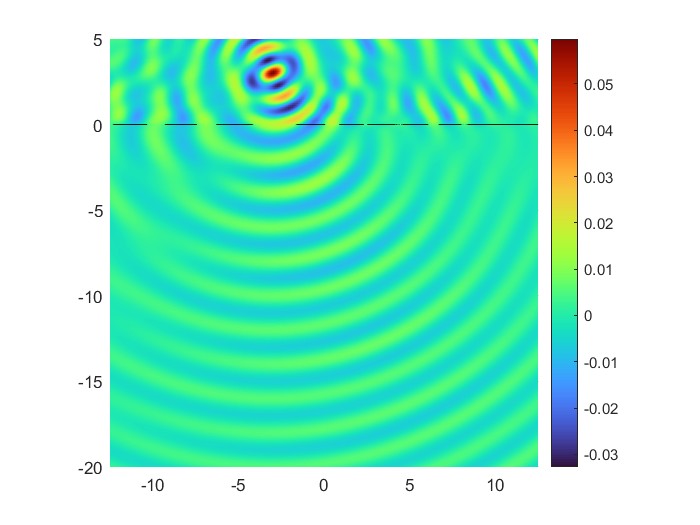}}
		\\
		\subfloat[]{
			\includegraphics[scale=0.25]{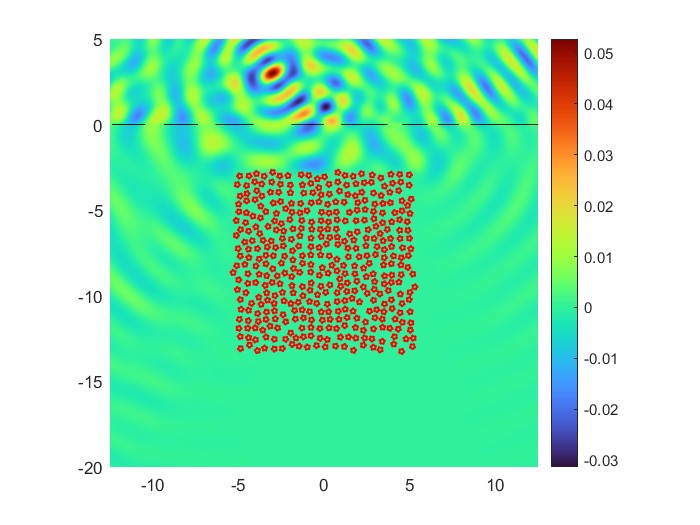}}
		\subfloat[]{
			\includegraphics[scale=0.25]{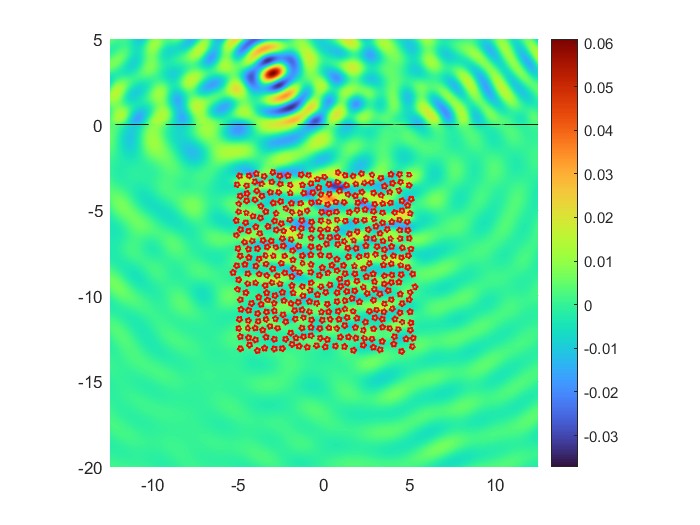}}
		\caption{Elastic scattering by 400 particles under point-source illumination with $\omega = 2\pi$. The imaginary part of the first component of the total elastic field is shown. (a) Geometry of the particles. (b) Incident point-source field without particles. (c) Field for multiple rigid particles. (d) Field for multiple traction-free particles.}
		\label{f4}
	\end{figure}

	In Figure \ref{f4}, we consider $M=400$ particles randomly distributed in $[-5,5]\times[-13,-3]$. The prototype particle is starfish-shaped, with parameters $a_1 = 0.135$, $a_2 =0.034$, and $a_3 =5$. See Figure \ref{f4}(a). Figure \ref{f4}(b) shows the incident point source wave in the two-layered background medium. The imaginary parts of the first component of the total fields are shown in Figure \ref{f4}(c) for the Dirichlet problem and Figure \ref{f4}(d) for the Neumann problem. We observe that when the average distance between particles is small, the scattered field behaves as if a large effective obstacle were present,
	particularly in the rigid case. 
	\begin{table}[h!]
		\begin{center}
			\caption{Results for elastic scattering by multiple buried starfish-shaped particles using the scattering matrix based method under point-source incidence.}
			\label{t3}
			\begin{tabular}{c|cc|ccc|ccc}
				\hline
				& & &\multicolumn{3}{c|}{Dirichlet}&\multicolumn{3}{c}{Neumann}\\
				\hline
				$\omega$&$N_{particle}$&$N_{tot}$&$N_{iter}$&$T_{solve}$&$E_{err}$&$N_{iter}$&$T_{solve}$&$E_{err}$\\
				\hline
				\multirow{3}{*}{$\pi$} 
				&100&2200&83&2.63E0&7.82E-10&136&5.10E0&1.26E-6\\
				&225&4950&87&7.40E0&4.35E-10&87&7.31E0&9.67E-7\\
				&400&8800&93&1.80E1&4.69E-10&75&1.26E1&5.78E-7\\
				\hline
				\multirow{3}{*}{$2\pi$} 
				&100&2200&175&6.97E0&1.93E-9&222&9.71E0&7.65E-7\\
				&225&4950&164&1.74E1&5.16E-10&460&9.31E1&8.51E-7\\
				&400&8800&179&4.09E1&5.79E-10&439&1.86E2&1.24E-6\\
				\hline
				\multirow{3}{*}{$4\pi$} 
				&100&2200&641&5.99E1&7.73E-8&420&2.97E1&1.59E-6\\
				&225&4950&703&1.92E2&2.97E-9&754&2.23E2&8.14E-7\\
				&400&8800&515&2.33E2&1.55E-9&829&5.73E2&8.73E-7\\
				\hline
			\end{tabular}
		\end{center}
	\end{table}

	\subsection{Scattering by plane wave incidence}
\begin{figure}[!t]
		\centering
		\subfloat[]{
			\includegraphics[scale=0.25]{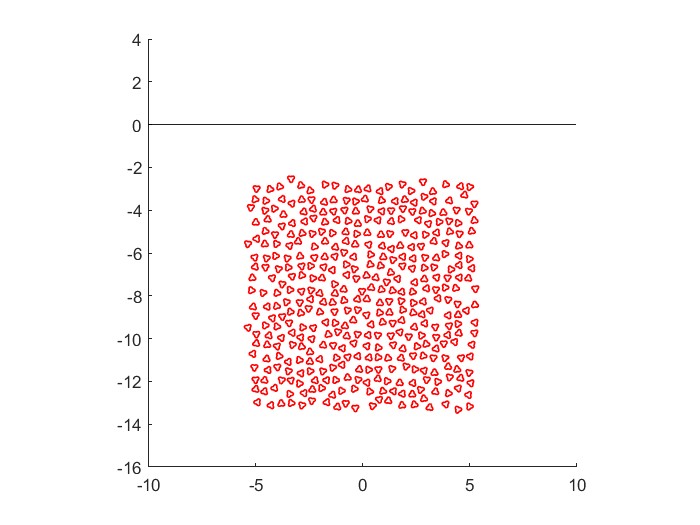}}
		\subfloat[]{
			\includegraphics[scale=0.25]{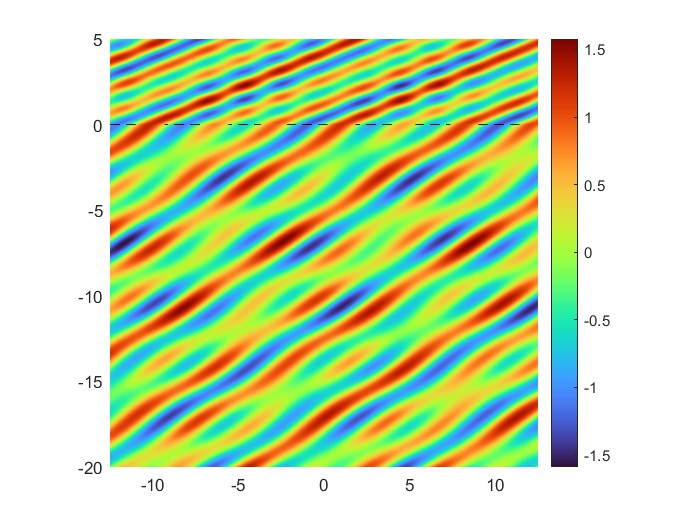}}
		\\
		\subfloat[]{
			\includegraphics[scale=0.25]{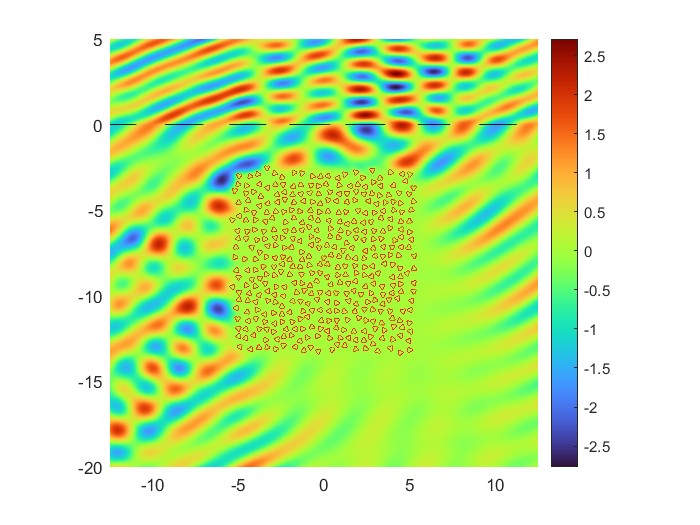}}
		\subfloat[]{
			\includegraphics[scale=0.25]{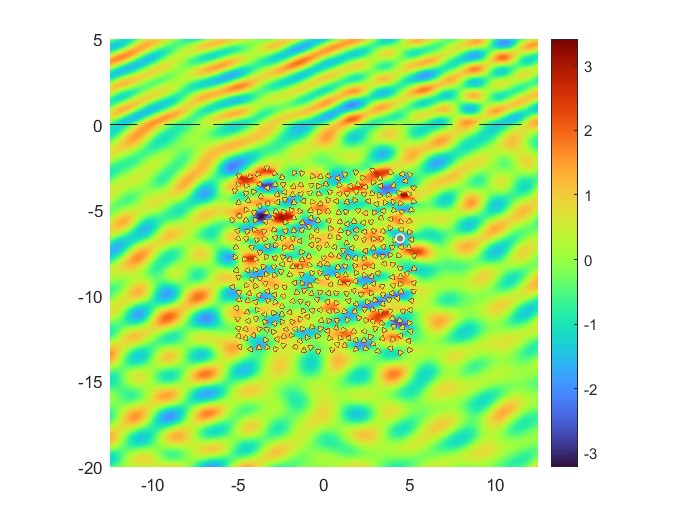}}
		\caption{Elastic scattering by 400 particles under plane-wave illumination with $\omega = 2\pi$. The imaginary part of the first component of the total elastic field is shown. (a) Geometry of the particles. (b) Plane-wave field without particles. (c) Field for multiple rigid particles. (d) Field for multiple traction-free particles.}
		\label{f5}
	\end{figure}

	In this example, we consider an incident plane wave defined by
	\[\uinc = \alpha\me^{\mi \kp\alpha\cdot x}+\alpha^\bot\me^{\mi \ks\alpha\cdot x},\quad \alpha = \left(\sin\left(\frac{\pi}{8}\right),-\cos\left(\frac{\pi}{8}\right)\right)^\top.\]
	We take $M=400$ particles randomly distributed in the square domain $[-5,5]\times[-13,-3]$. The prototype elastic particle is pear-shaped, with parameters $a_1 = 0.168$, $a_2 =0.042$, and $a_3 =3$. See Figure \ref{f5}(a) for the geometry. For plane-wave incidence, the right-hand side of the point-source linear system \eqref{aequ} must be modified. 
	
	First, we solve for the reflected field $\bu^r$ and transmitted field $\bu^t$ by imposing the penetrable boundary condition on the interface $\face$. Related derivations for Dirichlet and Neumann conditions are considered in \cite{brunoWindowedGreenFunction2021}. The background incident field $\bu^0$ in the two-layered elastic medium is then
	\begin{align}\label{bplan}
		\bu^0=\begin{cases}\uinc+\bu^r,&x\in\mathbb{R}_+^2,\\ \bu^t,&x\in\mathbb{R}_-^2.\end{cases}
	\end{align}
	It is classical that both $\bu^r$ and $\bu^t$ remain plane waves. We then construct the right-hand side of \eqref{aequ} by setting $\mathbf{b}=\mathbf{0}$ and replacing the lower zero block by the plane-wave expansion of $\bu^t$ on each buried particle.
	
	Table \ref{t4} presents numerical results for pear-shaped particles at various angular frequencies under plane wave incidence. Comparing Tables \ref{t3} and \ref{t4}, we find that the iteration counts for plane wave incidence are similar to those for point source incidence across different angular frequencies and boundary conditions. The results for both types of incidence also show that the GMRES iteration count depends on the particle size but is nearly independent of the particle shape. Figure \ref{f5}(b) shows the incident plane wave \eqref{bplan} in the two-layered background medium. The imaginary parts of the first component of the total fields are displayed in Figure \ref{f5}(c) for the Dirichlet problem and Figure \ref{f5}(d) for the Neumann problem. We use a larger plotting domain to show the wavefield behavior both inside and outside the particle cluster. In the rigid case, the incident field barely penetrates the cluster, whereas in the traction-free case the scattered field penetrates more deeply.
	\begin{table}[h!]
		\begin{center}
			\caption{Results for elastic scattering by multiple buried pear-shaped particles using the scattering-matrix-based method under plane-wave incidence.}
			\label{t4}
			\begin{tabular}{c|cc|ccc|ccc}
				\hline
				& & &\multicolumn{3}{c|}{Dirichlet}&\multicolumn{3}{c}{Neumann}\\
				\hline
				$\omega$&$N_{particle}$&$N_{tot}$&$N_{iter}$&$T_{solve}$&$E_{err}$&$N_{iter}$&$T_{solve}$&$E_{err}$\\
				\hline
				\multirow{3}{*}{$\pi$} 
				&100&2200&84&2.09E0&9.89E-9&103&2.98E0&3.75E-6\\
				&225&4950&87&7.28E0&3.69E-9&65&4.77E0&3.40E-6\\
				&400&8800&92&1.65E1&2.97E-9&57&9.50E0&3.17E-6\\
				\hline
				\multirow{3}{*}{$2\pi$} 
				&100&2200&256&1.17E1&2.47E-7&213&9.41E0&3.49E-6\\
				&225&4950&164&1.71E1&8.60E-9&392&6.90E1&6.15E-6\\
				&400&8800&176&4.23E1&1.29E-8&281&9.04E1&1.18E-5\\
				\hline
				\multirow{3}{*}{$4\pi$} 
				&100&2200&620&5.46E1&1.16E-6&345&2.03E1&1.09E-5\\
				&225&4950&918&3.12E2&1.22E-7&610&1.52E2&8.21E-6\\
				&400&8800&506&2.39E2&5.46E-8&815&5.56E2&9.34E-6 \\
				\hline
			\end{tabular}
		\end{center}
	\end{table}

    \subsection{Three-dimensional layered scattering}
\begin{figure}[!t]
		\centering
	\includegraphics[scale=0.8]{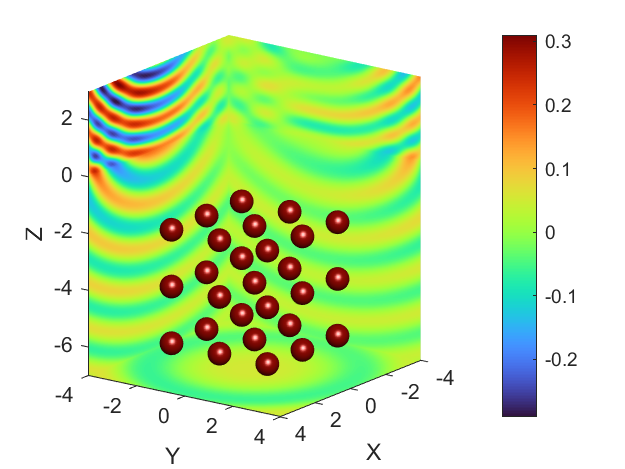}
		\caption{The incident point source wave and the distribution of 27 spheres.}
		\label{p1}
	\end{figure}

In this example, we present a 3D numerical experiment for layered elastic scattering by multiple buried particles.
In particular, we take $\omega = \frac{\pi}{2}$, $\lambda_1=0.2$, $\lambda_2=0.5$, $\mu_1 = 0.04$, and $\mu_2=0.1$, resulting in $\kp=2.97$, $\kp'=1.88$, $\ks=7.85$, and $\ks'=4.97$. The incident point source is located at $(0,0,6)$ with polarization direction $$\boldsymbol{p}=\left(\sin(\pi/3)\cos(\pi/4),\sin(\pi/3)\sin(\pi/4),\cos(\pi/3)\right).$$ 
We consider 27 rigid elastic spheres with Dirichlet boundary conditions in the lower half-space, each with radius $R=0.4$. The truncation number in the multipole expansion is chosen as $N=10$. For the numerical discretization of the 3D Sommerfeld integral, the contour parameters in \eqref{condd1} are chosen as $\alpha=1$, $\beta=1$, and $T_{max}=10$. The number of quadrature nodes is $N_r=48$ in the radial direction and $N_\theta=24$ in the angular direction.
The incident point source, corresponding to the 3D layered elastic Green's function, and the spheres are given in Figure \ref{p1}. In this experiment, all plotted fields show the real part of the first component.

\begin{figure}[!t]
		\centering
		\subfloat[]{
	\includegraphics[scale=0.4]{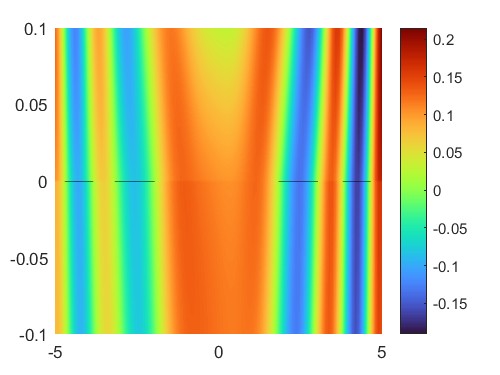}}
    \subfloat[]{
	\includegraphics[scale=0.4]{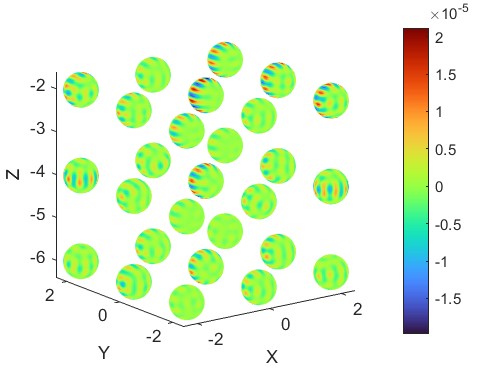}}
		\caption{(a) Total field across the interface. (b) Total field on the particles.}
		\label{p2}
	\end{figure}

\begin{figure}[!t]
		\centering
		\subfloat[]{
	\includegraphics[scale=0.45]{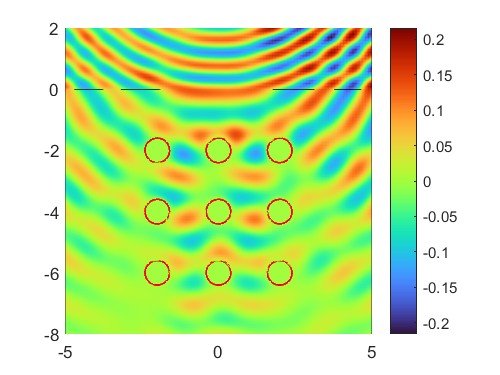}}
		\subfloat[]{
	\includegraphics[scale=0.45]{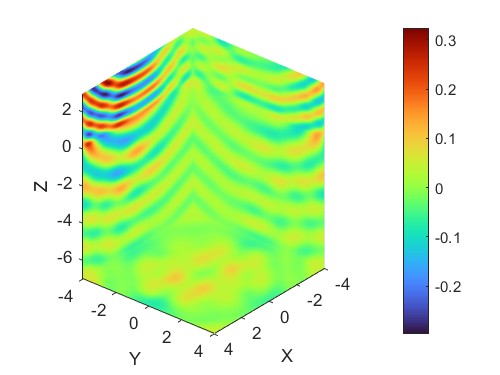}}\\
    \subfloat[]{
	\includegraphics[scale=0.45]{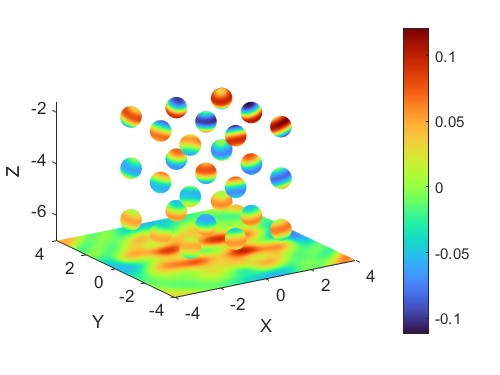}}
    \subfloat[]{
	\includegraphics[scale=0.45]{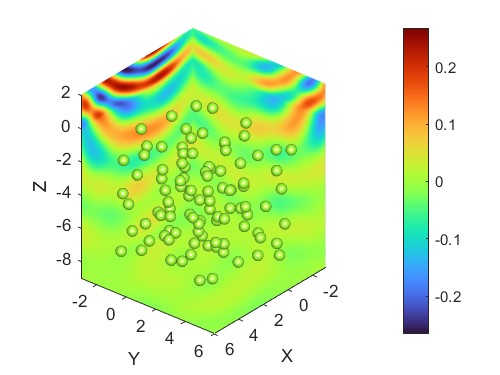}}
		\caption{(a) Total field on the slice $x=-2$. (b) Total field on three slices $x=-4$, $y=-4$, and $z=-7$. (c) Scattered field on the particles and on the slice $z=-7$. (d) Scattering by 100 particles with Dirichlet boundary conditions ($\omega={\pi}/{4}$).}
		\label{p3}
	\end{figure}

To validate the accuracy of the algorithm, we test the continuity of the computed total field across the interface by comparing the upward and downward total fields on the interface in Figure \ref{p2}(a). The results show that the field is continuous. We also plot the total field on the particle surfaces to verify the Dirichlet boundary condition in Figure \ref{p2}(b). The magnitude of the total field is of order $10^{-5}$, indicating good numerical accuracy.
Figure \ref{p3}(a) shows the total field on the slice $x=-2$, Figure \ref{p3}(b) shows the total field on three slices $x=-4$, $y=-4$, and $z=-7$, and Figure \ref{p3}(c) shows the scattered field on the particles and on the slice $z=-7$. Figure \ref{p3}(d) presents the total field for scattering by 100 buried particles with Dirichlet boundary conditions, where the computed field on the particles is again nearly zero. 

	\subsection{Inverse elastic scattering in  layered medium}
\begin{figure}
		\centering
		\includegraphics[scale =0.4]{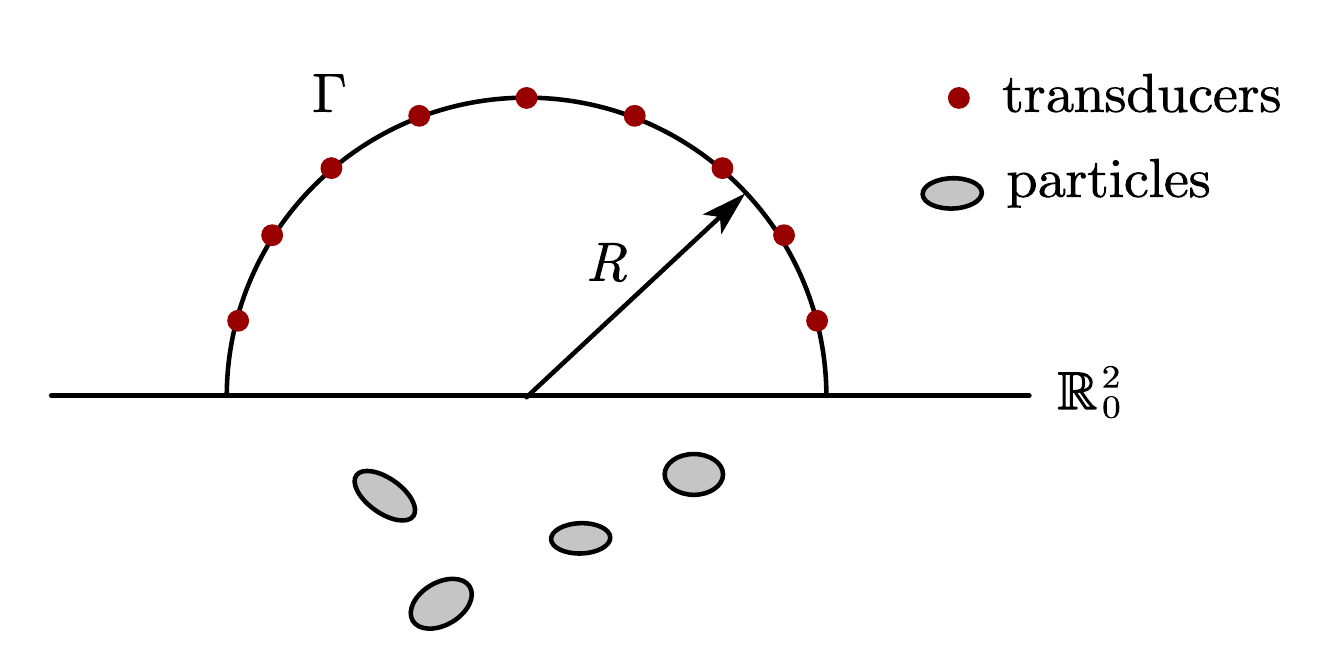}
		\caption{Geometry of the TRM.}\label{f3}
	\end{figure}

	The final numerical example illustrates an inverse scattering application of the proposed method. The goal is to localize the particles in the lower layer and generate selectively focused incident fields from the measurement of scattered field taken in the upper layer (see Figure \ref{f3}). To this end, we adopt a time-reversal technique known as the DORT method, which is based on the reversibility of the wave equation and has been rigorously justified in homogeneous elastic media with Dirichlet \cite{laiFastInverseElastic2022} and Neumann boundary conditions \cite{zhangSelectiveFocusingElastic}, as well as in layered media \cite{laizhang2026}. 

	In the near-field time-reversal process, the time-reversal mirror (TRM) $\Gamma$ emits a linear combination of point sources and measures the scattered field. More precisely, an incident wave $\uinc$ with density parameter $\boldsymbol{f}$ is given by
	\begin{align}\label{her}
		\uinc(\bx) = \int_\Gamma \Phi(\bx,\by)\boldsymbol{f}(\by)ds_{\by}
	\end{align}
	Let $\mathscr{F}\boldsymbol{f}$ denote the scattered field on $\Gamma$ generated by the incident wavefield with density $\boldsymbol{f}$. The time-reversal operator (TRO) is then given by $\mathscr{T} = \mathscr{F}^*\mathscr{F}$,
	where $\mathscr{F}^*$ is the adjoint operator of $\mathscr{F}$.
	DORT theory states that an eigenfunction $\boldsymbol{g}$ of $\mathscr{T}$ associated with a significant eigenvalue can be used as the density in \eqref{her} to generate an incident wave that selectively focuses on an unknown particle.
	
\begin{figure}[!t]
		\centering
		\subfloat[]{
			\includegraphics[scale=0.4]{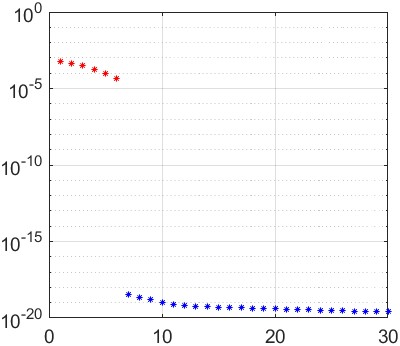}}
		\subfloat[]{
			\includegraphics[scale=0.4]{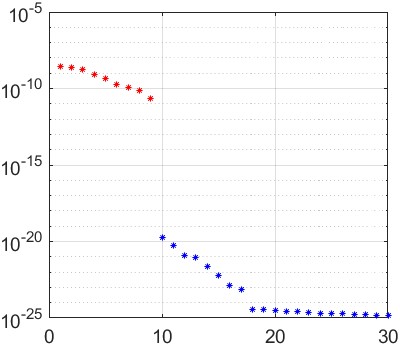}}
		\caption{Layered time reversal for elastic particles. (a) The first 30 largest eigenvalues of the time-reversal operator $T$ for rigid particles. (b) The first 30 largest eigenvalues of $T$ for traction-free particles.}
		\label{f7}
	\end{figure}

\begin{figure}[!t]
		\centering
	\subfloat[]{
		\includegraphics[scale=0.2]{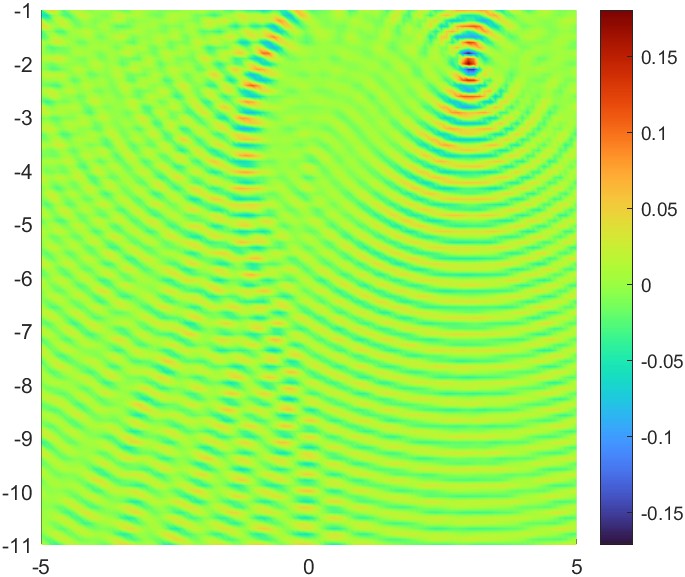}}
	\subfloat[]{
		\includegraphics[scale=0.2]{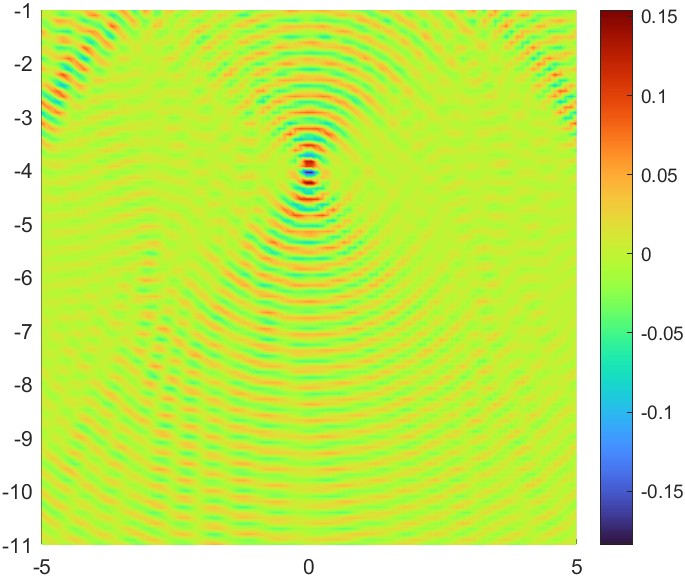}}
	\subfloat[]{
		\includegraphics[scale=0.2]{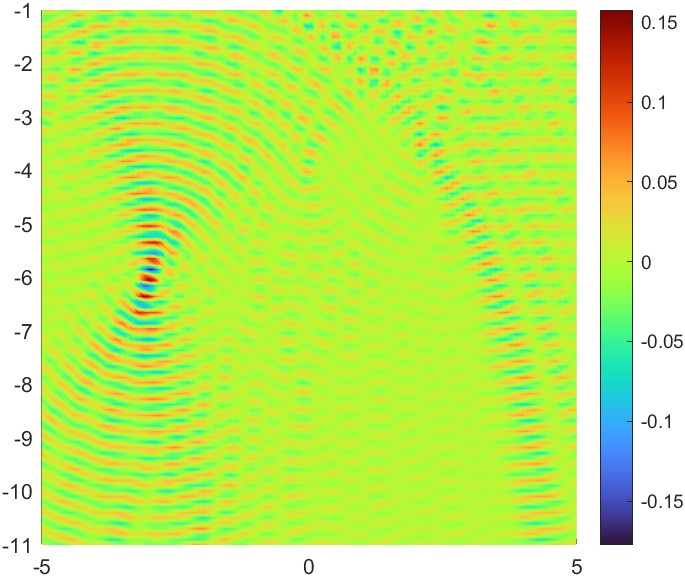}}\\
	\subfloat[]{
		\includegraphics[scale=0.2]{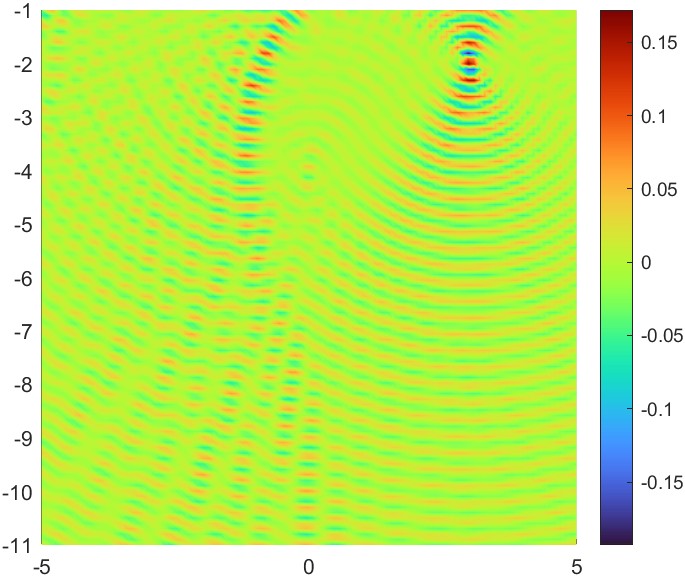}}
	\subfloat[]{
		\includegraphics[scale=0.2]{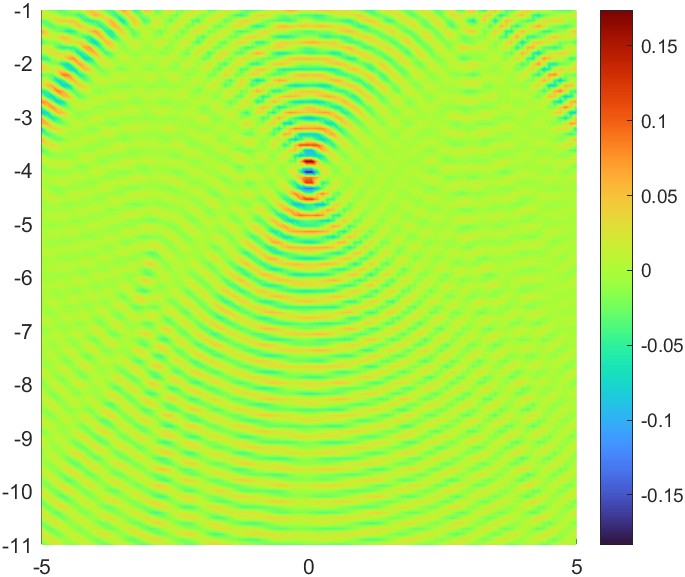}}
	\subfloat[]{
		\includegraphics[scale=0.2]{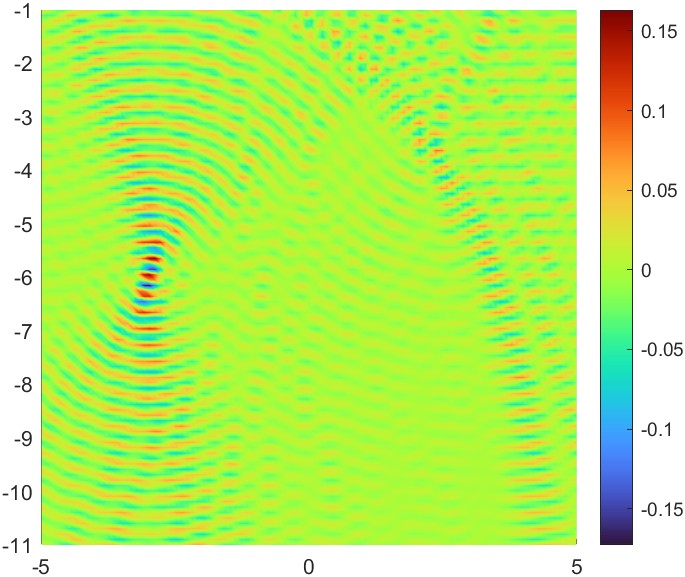}}
		\caption{Layered time reversal for elastic particles. The three (unknown) particles are located at $(3,-2)$, $(0,-4)$, and $(-3,-6)$. (a-c) Selective focusing of rigid particles. (d-f) Selective focusing of traction-free particles.}
		\label{f6}
	\end{figure}

	In the numerical experiment, we set the frequency to $\omega=40$ to achieve sufficient DORT resolution. We choose the TRM to be a semicircle $\Gamma$ with a radius $R=6$.
	The TRM is discretized using $N_\alpha = 30$ points, yielding a matrix approximation $F$ of $\mathscr{F}$ of size $2N_\alpha\times2N_\alpha$. The time-reversal matrix
	$T$ follows from the relation $T = \overline{F}^\top F$.
	We consider three rigid starfish-shaped particles located at $(3,-2)$, $(0,-4)$, and $(-3,-6)$, with sizes $3\times10^{-3}$, $2\times10^{-3}$, and $10^{-3}$, respectively. As shown in Figure \ref{f7}(a), the time-reversal matrix $T$ has six significant eigenvalues, with each rigid particle contributing two, consistent with results in homogeneous elastic media \cite{laiFastInverseElastic2022}. Figure \ref{f6}(a-c) shows the amplitudes of the incident wavefunctions \eqref{her} associated with the eigenfunctions corresponding to the first three eigenvalues. These wavefields selectively focus on the particles, although the focusing in the $y$ direction is slightly less sharp because of the limited aperture of the TRM in depth. Figure \ref{f6}(c) also shows that more deeply buried particles are reconstructed with lower resolution, a phenomenon also reported experimentally in acoustics \cite{devaneyTimeReversalImaging2005}. We repeat the same time-reversal procedure for traction-free particles, as shown in Figure \ref{f6}(d-f), and obtain similar focusing behavior. 
	
    
	\section{Conclusion}\label{conc}
	We have developed a fast algorithm for simulating elastic scattering by many rigid and traction-free particles buried in a layered medium under point-source and plane-wave incidence. The method combines Sommerfeld integral representations, which account for the layered interface, with a scattering matrix formalism for particles of arbitrary shape. By coupling these representations through Sommerfeld-to-local and multipole-to-Sommerfeld translation operators, the full layered multiple scattering problem is reduced to a structured linear system involving interface spectral densities and particle multipole coefficients. An effective Schur-complement preconditioner, together with fast multipole acceleration, allows the resulting system to be solved by GMRES with a modest number of iterations.

	Numerical experiments demonstrate the accuracy, efficiency, and flexibility of the proposed solver for both direct and inverse scattering problems. In particular, the method achieves high accuracy for noncircular particles, scales efficiently to large particle configurations, extends naturally to three-dimensional layered elastic scattering, and supports applications such as layered time-reversal focusing. Future work will focus on extending the framework to quasi-periodic structures and locally perturbed layered media.

	\bibliography{bib.bib}
    \bibliographystyle{abbrv} 
\end{document}